\documentclass[preprint,1p,leqno]{elsarticle}
\usepackage{amsmath}
\usepackage{amssymb}
\usepackage{mathrsfs}
\usepackage{bm}
\usepackage{color}
\usepackage[svgnames]{xcolor}
\usepackage{marvosym}
\usepackage{ulem}

\usepackage[pagewise]{lineno}

\usepackage[]{graphicx}
\usepackage{psfrag}
\usepackage{pdfpages}
\usepackage[]{xcolor}
\usepackage[makeroom]{cancel}

\usepackage{amsthm}

\usepackage{bigstrut}
\usepackage{multirow}
\usepackage{subfig}
\usepackage{easyReview}
\setreviewson

\usepackage{hyperref}

\usepackage{geometry}
\newcommand{\bu}{\boldsymbol{u}}
\newcommand{\bF}{\boldsymbol{f}}
\newcommand{\bw}{\boldsymbol{\omega}}

\newcommand{\bn}{\boldsymbol{n}}
\newcommand{\bj}{\boldsymbol{j}}

\newcommand{\bB}{\boldsymbol{B}}
\newcommand{\bH}{\boldsymbol{H}}
\newcommand{\bE}{\boldsymbol{E}}

\newcommand{\bv}{\boldsymbol{v}}

\newcommand{\Rn}{\mathrm{R}_f}
\newcommand{\Rm}{\mathrm{R}_m}
\newcommand{\Al}{\mathrm{A}_{l}}

\newcommand{\ip}[3]{\left\langle{#1},{#2}\right\rangle_{#3}}

\newcommand{\revAPalha}[1]{{\color{black}{#1}}}
\newcommand{\revABrugnoli}[1]{{\color{black}{#1}}}
\newcommand{\revYi}[1]{{\color{black}{#1}}}

\journal{arXiv}

\begin{document}
\begin{frontmatter}
\title{Decoupled structure-preserving discretization of incompressible MHD equations with general boundary conditions}

\author[1]{Yi Zhang\corref{cor1}}
\ead{zhangyi_aero@hotmail.com}
\author[2]{Artur Palha}
\ead{a.palha@tudelft.nl}
\author[3]{Andrea Brugnoli}
\ead{andrea.brugnoli@isae-supaero.fr}
\author[2]{Deepesh Toshniwal}
\ead{d.toshniwal@tudelft.nl}
\author[4]{Marc Gerritsma}
\ead{m.i.gerritsma@tudelft.nl}

\affiliation[1]{organization={Guilin University of Electronic Technology, School of Mathematics and Computing Science},
	city={Guilin},
	country={China}}

\affiliation[2]{organization={Delft Institute of Applied Mathematics, Delft University of Technology},
city={Delft},
country={the Netherlands}}

\affiliation[3]{organization={ICA, Université de Toulouse, ISAE–SUPAERO, INSA, CNRS, MINES ALBI, UPS},
	city={Toulouse},
	country={France}}
	
\affiliation[4]{organization={Faculty of Aerospace Engineering, Delft University of Technology},
city={Delft},
country={the Netherlands}}

\cortext[cor1]{Corresponding author}


\begin{abstract}
	In the framework of a mixed finite element method, a structure-preserving formulation for incompressible magnetohydrodynamic (MHD) equations with general boundary conditions is proposed. A leapfrog-type temporal scheme fully decouples the fluid part from the Maxwell part by means of staggered discrete time sequences and, in doing so, partially linearizes the system. Conservation and dissipation properties of the formulation before and after the decoupling are analyzed. We demonstrate optimal spatial and second-order temporal accuracy, as well as conservation and dissipation properties, of the proposed method using manufactured solutions, and apply it to the benchmark Orszag-Tang and lid-driven cavity cases.
\end{abstract}
\end{frontmatter}

\section{Introduction} \label{Sec: Introduction}
Given a bounded, contractible domain $ \Omega\subset\mathbb{R}^d $, $ d\in\left\lbrace2,3\right\rbrace $ with a Lipschitz boundary $ \partial\Omega $, the incompressible, constant density magnetohydrodynamic (or simply incompressible MHD) equations \cite{Davidson_2016} written as
\begin{subequations}\label{Eq: dimensional MHD}
	\begin{align}
		\label{dMHD a}
		\rho\left[ \partial_{t}\bu+\left(\bu\cdot\nabla \right) \bu\right]  - \tilde{\mu}\Delta\bu - \bj\times\bB + \nabla p &= \rho\bF , \\
		\label{dMHD f}\nabla\cdot \bu &=0 ,\\
		\label{dMHD b}\partial_{t}\bB + \nabla\times \bE  &= \boldsymbol{0} ,\\
		\label{dMHD c}\bj - \sigma\left( \bE +\bu\times\bB\right)  &= \boldsymbol{0} , \\
		\label{dMHD d}\bj - \nabla\times\boldsymbol{H} &=\boldsymbol{0} ,\\
		\label{dMHD e}\bB &= \mu\boldsymbol{H} ,
	\end{align}
\end{subequations}
govern the dynamics of fluid velocity $ \bu $, electric current density $ \bj $, magnetic flux density $ \bB $, magnetic field strength $ \boldsymbol{H} $, hydrodynamic pressure $ p $ and electric field strength $\bE$, subject to body force field $ \bF $, a velocity initial condition $\bu^{0} :=\bu(\boldsymbol{x}, 0)$ and a divergence-free magnetic initial condition $\bB^{0} :=\bB(\boldsymbol{x}, 0)$, and proper boundary conditions \cite{BENSALAH20015867}, in a space-time domain, $ \Omega\times(0,T] $. 
The material parameters involved are the fluid density $ \rho $, the dynamic viscosity $ \tilde{\mu} $,  the electric conductivity $ \sigma $, and the magnetic permeability $ \mu $. 

The incompressible MHD equations, \eqref{Eq: dimensional MHD}, are a coupled system of the incompressible Navier-Stokes equations and Maxwell's equations. The two evolution equations, \eqref{dMHD a} and \eqref{dMHD b}, stand for the conservation of momentum and Faraday’s law of induction, respectively. Ohm’s law and Ampere’s law are \eqref{dMHD c} and \eqref{dMHD d}, respectively. \eqref{dMHD e} sets up a constitutive relation between $ \bB$ and $ \bH $. The divergence-free constraint on $ \bu $, \eqref{dMHD f}, represents mass conservation of the constant mass density flow. It is seen that, as $ \bB^{0} :=\bB(\boldsymbol{x}, 0) $ is divergence-free, 
\eqref{dMHD b} ensures
\[
\nabla\cdot\bB = 0 \quad \mathrm{in}\ \Omega\times(0, T],
\]
which is Gauss’s law for magnetism. And, from \eqref{dMHD d}, we know $\nabla\cdot\bj = 0$ which stands for conservation of charge. The coupling between the flow field and the electromagnetic field is achieved through the terms $ \bj\times\bB  $ in \eqref{dMHD a} and $ \bu\times\bB $ in \eqref{dMHD c}.

Let
$ L $, $ U $, and $ B $
denote characteristic quantities of length, velocity, and magnetic flux density, respectively. We scale the space-time domain and variables by
\[
\boldsymbol{x} \leftarrow \dfrac{\boldsymbol{x}}{L},
\quad
t \leftarrow \dfrac{Ut}{L} ,
\]
\[
\bu \leftarrow \dfrac{\bu}{U},
\quad
p  \leftarrow \dfrac{p}{\rho U^2},
\quad
\bF \leftarrow \dfrac{ L\bF}{U^2} ,
\]
\begin{equation*}
		\bB  \leftarrow \dfrac{\bB}{B} ,
	\quad
	\bH  \leftarrow \dfrac{\mu\bH}{B} ,
	\quad
	\bE  \leftarrow \dfrac{\bE}{UB},
	\quad
	\bj  \leftarrow \dfrac{\mu L \bj}{B} .
\end{equation*}
And the induced scaled operators are
\[
D  \leftarrow L D, \  D\in\left\lbrace\nabla,\nabla\times,\nabla\cdot\right\rbrace,
\]
\[
\Delta \leftarrow L^2\Delta,
\quad 
\partial _{t} \leftarrow \dfrac{L}{U}\partial_{t},
\quad 
\bu \cdot\nabla \leftarrow \dfrac{L}{U}\left( \bu\cdot\nabla \right) .
\]
Note that we do not introduce new notation for the scaled objects, and, after this point, they always refer to the scaled dimensionless ones. 

Now \eqref{Eq: dimensional MHD} can be written in a dimensionless form, 
\begin{subequations}\label{Eq: dimensionless ast 0}
	\begin{align}
		\partial _{t}\bu +\bw\times \bu  + \Rn^{-1}\nabla\times \bw  - \mathsf{c}\left(\nabla\times\bH\right) \times\bH  + \nabla  P  &= \bF, &&\qquad\text{in }\Omega\times(0,T], \label{Eq: dimensionless ast 0 a}\\
		\bw - \nabla\times\bu &=\boldsymbol{0}, &&\qquad\text{in } \Omega\times(0, T],\label{Eq: dimensionless ast 0 b}\\
		\nabla\cdot   \bu  &=0,&&\qquad\text{in }\Omega\times(0,T],\label{Eq: dimensionless ast 0 c}\\
		\partial _{t}\bH  +  \Rm^{-1}\Delta  \bH  - \nabla\times\left(\bu \times\bH\right)  &= \boldsymbol{0},&&\qquad\text{in }\Omega\times(0,T], \label{eq: maxwell original}
	\end{align}
\end{subequations}
where, after introducing vorticity $\bw:=\nabla\times\boldsymbol{u}$ and total pressure $P:= p + \dfrac{1}{2}\bu\cdot\bu$, we have applied the relation
$
\left( \bu\cdot\nabla\right) \bu + \nabla p\stackrel{\eqref{Eq: dimensionless ast 0 c}}{=} \bw\times \bu + \nabla P
$.
Meanwhile, by eliminating $\bB$, $\bE$ and $\bj$, dimensionless forms of \eqref{dMHD b}-\eqref{dMHD e}, i.e.
\begin{subequations}\label{eq: compress}
	\begin{align}
		\partial_{t}\bB + \nabla\times \bE  &= \boldsymbol{0} ,\label{eq: BE}\\
		\Rm^{-1}\bj - \left( \bE +\bu\times\bB\right)  &= \boldsymbol{0} , \\
		\bj - \nabla\times\boldsymbol{H} &=\boldsymbol{0} ,\label{eq: compress c}\\
		\bB &= \bH,\label{eq: compress d}
	\end{align}
\end{subequations}
have been written into \eqref{eq: maxwell original},
which is the magnetic induction (or Helmholtz) equation.\footnote{Note that, after non-dimensionalization, we get $\bB=\bH$, i.e. \eqref{eq: compress d}. This is the case in vector calculus. In differential forms, they appear as a pair of differently oriented dual forms connected by the Hodge $\star$ operator, i.e. $\bB = \star\bH$.}
The three dimensionless parameters are the fluid Reynolds number $ \Rn=\dfrac{\rho U L}{\tilde{\mu}} = \dfrac{U L}{\nu} $ (with $ \nu=\dfrac{\tilde{\mu}}{\rho} $ being the kinematic viscosity), the magnetic Reynolds number $ \Rm=\mu\sigma U L $, and the coupling number $\mathsf{c} = \Al^{-2}$ ($ \Al  = \dfrac{U\sqrt{\rho\mu}}{B}= \dfrac{U}{U_{\mathrm{A}}} $  is the Alfvén number and $ U_{\mathrm{A}} = \dfrac{B}{\sqrt{\rho\mu}} $ is the Alfvén speed). 
Sometimes, the coupling term in \eqref{Eq: dimensionless ast 0 a}, $\mathsf{c}\left(\nabla\times\bH \right)  \times\bH$, is written as 
$\mathsf{s}\left( \bE  +\bu \times\bH \right) \times\bH$
where the coupling number becomes 
$\mathsf{s}=\mathsf{c}\Rm = \dfrac{\sigma L B^2}{\rho U}$. 
For example, see \cite[(2.2a)]{ZHANG202345}. 
\revAPalha{One can see that, from \eqref{eq: compress}, 
\[
    \partial_{t}\boldsymbol{H} = \partial_{t}\boldsymbol{B} = -\nabla\times \boldsymbol{E},
\]
and due to this, we will have that
\begin{equation}\label{eq: divH}
    \partial_{t}\nabla\cdot\boldsymbol{H} = -\nabla\cdot\nabla\times \boldsymbol{E} = 0,
\end{equation}
which leads to Gauss's law for magnetism if $\boldsymbol{H}^0$ is divergence-free.
}

To close the problem, \eqref{Eq: dimensionless ast 0} is complemented with initial conditions
$\bu^{0}, \bH^{0}$
and boundary conditions. In this work, we consider the following general boundary condition configuration:
\begin{itemize}
	\item For the \emph{fluid} part, there are two pairs of boundary conditions,
	\begin{equation}\label{eq: bc NS}
		\left\lbrace
		\begin{aligned}
			&P = \widehat{P} \quad&&\text{on }\Gamma_{\widehat{P}}\times\left(0,T\right]\\
			&\bu\cdot\boldsymbol{n} = \hat{u}\quad&&\text{on }\Gamma_{\hat{u}}\times\left(0,T\right]
		\end{aligned}\right.
		\quad \text{and}\quad 
		\left\lbrace
		\begin{aligned}
			&\bu\times\boldsymbol{n} = \widehat{\bu}\quad&&\text{on }\Gamma_{\widehat{\bu}}\times\left(0,T\right]\\
			&\bn\times\left( \boldsymbol{\omega}\times\bn\right)  = \widehat{\bw} \quad&&\text{on }\Gamma_{\widehat{\bw}}\times\left(0,T\right]
		\end{aligned}\right.,
	\end{equation}
	where $ \bn$ is the outward unit normal vector, $\left\lbrace\Gamma_{\widehat{P}}, \Gamma_{\hat{u}}\right\rbrace$ and $\left\lbrace\Gamma_{\widehat{\bu}}, \Gamma_{\widehat{\bw}}\right\rbrace$ are two partitions of the boundary, such that $\overline{\Gamma}_{\widehat{P}}\cup \overline{\Gamma}_{\hat{u}} = \partial\Omega$, $\Gamma_{\widehat{P}}\cap \Gamma_{\hat{u}} = \emptyset$ and $\overline{\Gamma}_{\widehat{\bu}}\cup \overline{\Gamma}_{\widehat{\bw}} = \partial\Omega$, $\Gamma_{\widehat{\bu}}\cap \Gamma_{\widehat{\bw}} = \emptyset$.
	
	\item For the \emph{Maxwell} part, boundary conditions are given as
	\begin{equation}\label{eq: bc Maxwell}
		\left\lbrace
		\begin{aligned}
			&\boldsymbol{E}\times\bn = \widehat{\boldsymbol{E}} \quad&&\text{on }\Gamma_{\widehat{\boldsymbol{E}}}\times\left(0,T\right]\\
			&\bn\times\left( \bH\times\boldsymbol{n}\right)  = \widehat{\bH}\quad&&\text{on }\Gamma_{\widehat{\bH}}\times\left(0,T\right]
		\end{aligned}\right.,
	\end{equation}
	where $\left\lbrace\Gamma_{\widehat{\boldsymbol{E}}}, \Gamma_{\widehat{\bH}}\right\rbrace$ is a partition of the boundary, satisfying $\overline{\Gamma}_{\widehat{\boldsymbol{E}}}\cup \overline{\Gamma}_{\widehat{\bH}} = \partial\Omega$, $\Gamma_{\widehat{\boldsymbol{E}}}\cap \Gamma_{\widehat{\bH}} = \emptyset$.
	Note that the boundary condition $\boldsymbol{E}\times\bn $ implies the magnetic flux boundary condition $\bH\cdot\bn$ (or $\bB\cdot\bn$). This is because
\begin{equation}\partial_{t}\left( \bH\cdot\boldsymbol{n}\right) =\partial_{t}\left( \bB\cdot\boldsymbol{n}\right) \stackrel{\eqref{eq: BE}}{=} - \left( \nabla\times \bE \right) \cdot \boldsymbol{n} \stackrel{\nabla\times\bn=\boldsymbol{0}}{=} -\left( \nabla\cdot\right) _{\partial\Omega}\left( \boldsymbol{E}\times \boldsymbol{n}\right) \label{eq:H_perp_E_tan}\end{equation}
on $\partial\Omega\times\left(0,T\right]$. Also, see \cite[Section~2.3]{hu2017divB}. 

\end{itemize}
For example, the boundary conditions considered in the work of Hu et al., see (7) of \cite{HU2021110284}, are a particular case of this general configuration, i.e. the case when $\Gamma_{\hat{u}} = \Gamma_{\widehat{\bw}} = \Gamma_{\widehat{\bH}}=\emptyset$,  $\widehat{P} = 0 $ and $\widehat{\bu}=\widehat{\boldsymbol{E}}=\boldsymbol{0}$.

It is generally acknowledged that structure-preservation, especially the preservation (either strongly or weakly) of $\nabla\cdot\bB =0$, is an essential property of numerical methods for MHD \cite{BRACKBILL1980426}. Many of these structure-preserving methods can be classified into, for example, potential-based methods \cite{DING2024}, divergence-cleaning methods \cite{TOTH2000605,Balsara_2004}, constrained transport methods \cite{Evans1988SimulationOM,refId0} and Lagrange multiplier (or augmented) methods \cite{Heister2017,ZHANG2022110752,10.1093/imanum/drad005}. In the current paper, none of those techniques will be applied, but it will be shown that the absence of magnetic monopoles is weakly enforced through the decoupled variational formulation.

Recently, methods that aim at preserving multiple structures in the mixed finite element setting have become increasingly popular. See, for example, the work of Hu et al. for incompressible MHD that preserves cross- and magnetic-helicity, energy, and Gauss law of magnetism \cite{HU2021110284}, the work of Gawlik and Gay-Balmaz for incompressible MHD of a variable fluid density that preserves energy, cross-helicity (when the fluid density is constant) and magnetic-helicity, mass, total squared density, pointwise incompressibility, and Gauss's law for magnetism \cite{GAWLIK2022110847}, the work of Laakmann et al. \cite{LAAKMANN2023112410} and references therein. These methods preserve several physical quantities of interest, but are usually less computationally efficient because of their large (nonlinear) discrete systems to be solved. 

A compromise is to use a decoupled scheme that breaks the discrete system into several smaller ones. \revYi{For example, Zhang and Su \cite{ZHANG202345} propose a three-step decoupled scheme by using a first-order semi-implicit Euler scheme with a first-order stabilizer and some carefully designed implicit-explicit treatments for the coupling terms. In \cite{ZHANG2022110752}, a four-step decoupled scheme of second-order temporal accuracy is studied. It employs the mixed finite element method under a \textquotedblleft zero-energy-contribution\textquotedblright\ idea and a projection method. Ma and Huang introduced a vector penalty-projection method, a five-step decoupled scheme that has a mixed finite element spatial discretization and a first-order backward Euler temporal discretization in \cite{MA202228}. The work \cite{ZHANG2024210} reports another five-step decoupled scheme which combines the idea of an auxiliary variable, pressure correction projection method, and some other implicit-explicit techniques for the thermally coupled MHD equations with variable density. For a more extended literature study on decoupled schemes in a mixed finite element setting for MHD equations, we refer to \cite{MA202228}.} These methods commonly deal with a particular configuration of boundary conditions.

In this paper, \revYi{we introduce a mixed finite element structure-preserving decoupled discretization of a second-order temporal accuracy for the incompressible MHD equations with general boundary conditions. This method decouples each full iteration into two steps by means of an implicit mid-point scheme that is applied to the two evolution equations at staggered temporal sequences such that they can be solved individually in a leapfrog approach.}

The outline of the rest of the paper is as follows. Relevant Hilbert spaces and the de Rham Hilbert complex are introduced in Section~\ref{Sec: Preliminaries} followed by the proposed semi-discrete formulation and its conservation properties in Section~\ref{Sec: formulation}. The temporal discretization that leads to the decoupled scheme is explained in Section~\ref{Sec: discretization}. Numerical tests are presented in Section~\ref{Sec: numerical tests}. Finally, in Section~\ref{Sec: conclusion}, conclusions are drawn.



\section{Relevant Hilbert spaces and the de Rham complex} \label{Sec: Preliminaries}
The function spaces used in this work are similar to those used in \cite{ZHANG2024113080} where only the two-dimensional setting is discussed. For the sake of completeness, we give a brief introduction to the three-dimensional spaces as a complement to \cite[Section~2.1]{ZHANG2024113080}.
We first introduce the relevant spaces and the de Rham complex in the infinite-dimensional setting and then introduce their finite-dimensional counterparts.
For an extensive introduction to function spaces, see \cite{Oden2010} or any other textbook on functional analysis.

\subsection{The infinite-dimensional setting}
Let $L^2\left(\Omega\right)$ denote the space of square-integrable functions on $\Omega$,
\[
L^2(\Omega) := \left\lbrace\varphi\left|\ip{\varphi}{\varphi}{\Omega} :=\int_{\Omega}\varphi\cdot\varphi\ \mathrm{d}\Omega<\infty\right.\right\rbrace.
\]
For $\Omega \subset \mathbb{R}^3$, we also consider the following Hilbert spaces,
\[
\begin{aligned}
H^1\left(\Omega\right) &:= \left\lbrace\psi\left|\psi\in L^2(\Omega),\nabla\psi\in\left[L^2\left(\Omega\right)\right]^3 \right.\right\rbrace,
\\
H\left(\mathrm{curl};\Omega\right) &:= \left\lbrace\boldsymbol{\sigma}\left|\boldsymbol{\sigma}\in \left[L^2\left(\Omega\right)\right]^3,\nabla\times\boldsymbol{\sigma}\in\left[L^2\left(\Omega\right)\right]^3 \right.\right\rbrace,
\\
H\left(\mathrm{div};\Omega\right) &:= \left\lbrace\boldsymbol{\phi}\left|\boldsymbol{\phi}\in \left[L^2\left(\Omega\right)\right]^3,\nabla\cdot\boldsymbol{\phi}\in L^2\left(\Omega\right) \right.\right\rbrace.
\end{aligned}
\]
They form the well-known de Rham Hilbert complex \cite{Arnold_Falk_Winther_2006, boffi2013mixed},
\[
	\mathbb{R}\hookrightarrow H^1\left(\Omega\right) \stackrel{\nabla}{\longrightarrow} H(\mathrm{curl};\Omega) \stackrel{\nabla\times}{\longrightarrow}  H(\mathrm{div};\Omega) \stackrel{\nabla\cdot}{\longrightarrow}  L^2(\Omega) \to 0.
\]
The above complex is exact in the contractible domain $Omega$, i.e. the image of $\nabla$ coincides with the kernel of $\nabla \times$, and the image of $\nabla \times$ coincides with the kernel of $\nabla \cdot$.

The trace operator, $\mathcal{T}$, restricts an element of Hilbert spaces to a boundary section $\Gamma\subseteq\partial\Omega$. For $\psi\in H^1\left(\Omega\right)$ and $\boldsymbol{\phi}\in H\left(\mathrm{div};\Omega\right) $, we define its action as
\[\mathcal{T}\psi = \left.\psi\right|_{\Gamma},\quad \mathcal{T}\boldsymbol{\phi} =\left.\boldsymbol{\phi}\cdot\bn\right|_{\Gamma}. \]
And we distinguish between $\mathcal{T}$ and $\mathcal{T}_{\parallel}$ for elements of $H\left(\mathrm{curl};\Omega\right) $, i.e., for $\boldsymbol{\sigma},\boldsymbol{\eta}\in H\left(\mathrm{curl};\Omega\right)$,
\begin{equation}\label{eq: curl traces}
	\mathcal{T}\boldsymbol{\sigma} = \left.\boldsymbol{\sigma}\times \bn\right|_{\Gamma},\quad \mathcal{T}_{\parallel}\boldsymbol{\eta} = \left.\bn\times \left( \boldsymbol{\eta}\times \bn\right)\right|_{\Gamma} .
\end{equation}
Using the above, we define the trace spaces on $\Gamma$ as
\[
\begin{aligned}
	H^{1/2}\left(\Omega,\Gamma\right):=\left\lbrace\mathcal{T}\psi\left|\psi\in H^1(\Omega)\right.\right\rbrace,
\end{aligned}
\]
\[
\begin{aligned}
	\mathcal{T}H\left(\mathrm{curl};\Omega,\Gamma\right):=\left\lbrace\mathcal{T}\boldsymbol{\sigma}\left|\boldsymbol{\sigma}\in H(\mathrm{curl};\Omega)\right.\right\rbrace,
\end{aligned}
\]
\[
\begin{aligned}
	\mathcal{T}_{\parallel}H\left(\mathrm{curl};\Omega,\Gamma\right):=\left\lbrace\mathcal{T}_{\parallel}\boldsymbol{\eta}\left|\boldsymbol{\eta}\in H(\mathrm{curl};\Omega)\right.\right\rbrace,
\end{aligned}
\]
\[
\begin{aligned}
	\mathcal{T}H\left(\mathrm{div};\Omega,\Gamma\right):=\left\lbrace\mathcal{T}\boldsymbol{\phi}\left|\boldsymbol{\phi}\in H(\mathrm{div};\Omega)\right.\right\rbrace.
\end{aligned}
\]
Note that, at the continuous level, $\mathcal{T}$ and $\mathcal{T}_{\parallel}$ are different, see \eqref{eq: curl traces}, and trace spaces $\mathcal{T}H\left(\mathrm{curl};\Omega,\Gamma\right)$ and $\mathcal{T}_{\parallel}H\left(\mathrm{curl};\Omega,\Gamma\right)$ \revABrugnoli{are in a duality, see \cite{buffa2001traces}.}

\subsection{The finite-dimensional setting}\label{subseq:finite_dim_setting}
Assume that we are given finite-dimensional spaces
$
G(\Omega)\subset  H^1(\Omega)$, $ C(\Omega)\subset  H(\mathrm{curl};\Omega)$, $D(\Omega)\subset H(\mathrm{div};\Omega)$, $S(\Omega)\subset L^2(\Omega)
$ which form a discretization of the continuous de Rham Hilbert complex,
\begin{equation} \label{eq: discrete de Rham complex}
	\mathbb{R}\hookrightarrow G(\Omega) \stackrel{\nabla}{\longrightarrow} C(\Omega) \stackrel{\nabla\times}{\longrightarrow}  D(\Omega) \stackrel{\nabla\cdot}{\longrightarrow} S(\Omega) \to 0,
\end{equation}
which is also exact in the contractible domain $Omega$.
We will also consider the following subspaces with boundary conditions,
where $\Gamma \subset \partial \Omega$, 
\[
G_{\widehat{\psi}}(\Omega,\Gamma):=\left\lbrace\psi_{h}\left|\psi_{h}\in G(\Omega),\ \mathcal{T}\psi_{h}=\widehat{\psi}\ \text{on}\ \Gamma\right.\right\rbrace,
\]
\[
C_{\boldsymbol{\widehat{\sigma}}}(\Omega,\Gamma):=\left\lbrace\boldsymbol{\sigma}_{h}\left|\boldsymbol{\sigma}_{h}\in C(\Omega),\ \mathcal{T}\boldsymbol{\sigma}_{h}=\boldsymbol{\widehat{\sigma}}\ \text{on}\ \Gamma\right.\right\rbrace,
\]
\[
C^{\parallel}_{\boldsymbol{\widehat{\eta}}}(\Omega,\Gamma):=\left\lbrace\boldsymbol{\eta}_{h}\left|\boldsymbol{\eta}_{h}\in C(\Omega),\ \mathcal{T}_{\parallel}\boldsymbol{\eta}_{h}=\boldsymbol{\widehat{\eta}}\ \text{on}\ \Gamma\right.\right\rbrace,
\]
\[
D_{\widehat{\phi}}(\Omega,\Gamma):=\left\lbrace\boldsymbol{\phi}_{h}\left|\boldsymbol{\phi}_{h}\in D(\Omega),\ \mathcal{T}\boldsymbol{\phi}_{h}=\widehat{\phi}\ \text{on}\ \Gamma\right.\right\rbrace.
\]
For the corresponding homogeneous boundary conditions we have
\[
G_{0}(\Omega,\Gamma):=\left\lbrace\psi_{h}\left|\psi_{h}\in G(\Omega),\ \mathcal{T}\psi_{h}=0\ \text{on}\ \Gamma\right.\right\rbrace,
\]
\[
C_{\boldsymbol{0}}(\Omega,\Gamma):=\left\lbrace\boldsymbol{\sigma}_{h}\left|\boldsymbol{\sigma}_{h}\in C(\Omega),\ \mathcal{T}\boldsymbol{\sigma}_{h}=\boldsymbol{0}\ \text{on}\ \Gamma\right.\right\rbrace,
\]
\[
C^{\parallel}_{\boldsymbol{0}}(\Omega,\Gamma):=\left\lbrace\boldsymbol{\eta}_{h}\left|\boldsymbol{\eta}_{h}\in C(\Omega),\ \mathcal{T}_{\parallel}\boldsymbol{\eta}_{h}=\boldsymbol{0}\ \text{on}\ \Gamma\right.\right\rbrace,
\]
\[
D_{0}(\Omega,\Gamma):=\left\lbrace\boldsymbol{\phi}_{h}\left|\boldsymbol{\phi}_{h}\in D(\Omega),\ \mathcal{T}\boldsymbol{\phi}_{h}=0\ \text{on}\ \Gamma\right.\right\rbrace.
\]
Moreover, the finite-dimensional trace spaces will be denoted by $\mathcal{T}G(\Omega,\Gamma)$, $\mathcal{T}C(\Omega,\Gamma)$, $\mathcal{T}_{\parallel}C(\Omega,\Gamma)$, $\mathcal{T}D(\Omega,\Gamma)$. Unlike their continuous counterparts, finite-dimensional trace spaces $\mathcal{T}C(\Omega,\Gamma)$ and $\mathcal{T}_{\parallel}C(\Omega,\Gamma)$ usually are different.

We define a corresponding trilinear form,
\begin{equation}\label{eq: trilinear}
	a\left(\boldsymbol{\alpha}_h,\boldsymbol{\beta}_h, \boldsymbol{\gamma}_h\right) := \ip{\boldsymbol{\alpha}_h\times\boldsymbol{\beta}_h}{\boldsymbol{\gamma}_h}{\Omega}, \quad \boldsymbol{\alpha}_h, \boldsymbol{\beta}_h,\boldsymbol{\gamma}_h \in \left\lbrace  C(\Omega),  D(\Omega)\right\rbrace,
\end{equation}
which is skew-symmetric with respect to any two of the three entries. Note that, 
in the infinite-dimensional setting, the $L^2$-boundedness of $\ip{\boldsymbol{\alpha}\times\boldsymbol{\beta}}{\boldsymbol{\gamma}}{\Omega}$ is not implied even if $\boldsymbol{\alpha}$, $\boldsymbol{\beta}$ and $\boldsymbol{\gamma}$ are $L^2$ vector fields.
However, in the finite-dimensional setting where the spaces $C(\Omega)$ and $D(\Omega)$ consist of piecewise-polynomial finite element vector fields, it is in fact $L^2$-bounded.


\section{Spatially-discrete formulation and its conservation properties} \label{Sec: formulation}
In this section, we present a spatially discrete weak formulation of \eqref{Eq: dimensionless ast 0} and analyze its conservation properties.
As a remark, for the sake of neatness, we discuss only the three-dimensional formulation and its temporal discretization in this paper. Obtaining their two-dimensional versions is straightforward.

\subsection{Spatially-discrete formulation}
To increase the notational clarity, we omit the part for the temporal domain. The proposed semi-discrete weak formulation of \eqref{Eq: dimensionless ast 0} in $\mathbb{R}^3$ is written as: Given $\bF\in\left[L^2(\Omega)\right]^3$, natural boundary conditions $\widehat{P}\in H^{1/2}\left(\Omega;\Gamma_{\widehat{P}}\right),\  \widehat{\bu}\in\mathcal{T}H\left(\mathrm{curl};\Omega,\Gamma_{\widehat{\bu}}\right) ,\  \widehat{\bE}\in \mathcal{T}H\left(\mathrm{curl};\Omega,\Gamma_{\widehat{\boldsymbol{E}}}\right)$,
and initial conditions 
$\left(\bu_{h}^{0},\bH^{0}_{h}\right)\in D\left(\Omega\right)\times C\left(\Omega\right)$,
seek $\left(\boldsymbol{u}_{h}, \bw_{h}, P_{h}, \bH_{h}\right)\in D_{\hat{u}}(\Omega,\Gamma_{\hat{u}})\times C^{\parallel}_{\widehat{\bw}}(\Omega,\Gamma_{\widehat{\bw}})\times S(\Omega)\times C^{\parallel}_{\widehat{\boldsymbol{H}}}(\Omega,\Gamma_{\widehat{\boldsymbol{H}}})$, such that $\forall \left(\bv_{h},\boldsymbol{w}_{h}, q_{h},\boldsymbol{b}_{h}\right)\in D_{0}\left(\Omega,\Gamma_{\hat{u}}\right) \times C^{\parallel}_{\boldsymbol{0}}\left(\Omega,\Gamma_{\widehat{\bw}}\right)\times S\left(\Omega\right)\times C_{\boldsymbol{0}}^{\parallel}\left(\Omega,\Gamma_{\widehat{\bH}}\right) $,
\begin{subequations}\label{Eq: wf}
	\begin{align}
		\label{wf a}
		\left\langle\partial _{t}\bu_{h}, \bv_{h}\right\rangle_{\Omega}
		+
		a\left(\bw_{h}, \bu_{h}, \bv_{h}\right)
		 + \Rn^{-1}\ip{\nabla\times\bw_{h}}{\bv_{h}}{\Omega} \qquad &\\ 
		\nonumber
		- \mathsf{c}\ 
		a\left(\nabla\times\bH_{h},\bH_{h},\bv_{h}\right)
		- \ip{P_{h}}{\nabla\cdot\bv_{h}}{\Omega}  &= \ip{\bF}{\bv_{h}}{\Omega} - \ip{\widehat{P}}{\mathcal{T}\bv_{h}}{\Gamma_{\widehat{P}}},\\
		\label{wf b}
		-\ip{\bu_{h}}{\nabla\times\boldsymbol{w}_{h}}{\Omega}+\ip{\bw_{h}}{\boldsymbol{w}_{h}}{\Omega} & = - \ip{\widehat{\bu}}{\mathcal{T}_{\parallel}\boldsymbol{w}_{h}}{\Gamma_{\widehat{\bu}}} ,\\
		\label{wf c}
		\ip{\nabla\cdot\bu_{h}}{q_{h}}{\Omega} & = 0,\\
		\label{wf d}
		\ip{\partial_{t}\bH_{h}}{\boldsymbol{b}_{h}}{\Omega} + \Rm^{-1}\ip{\nabla\times \bH_{h}}{\nabla\times\boldsymbol{b}_{h}}{\Omega} 
		-\ a\left(\bu_{h}, \bH_{h}, \nabla\times\boldsymbol{b}_{h}\right)
		&= \ip{\widehat{\bE}}{\mathcal{T}_{\parallel}\boldsymbol{b}_{h}}{\Gamma_{\widehat{\bE}}}.
	\end{align}
\end{subequations}
Note that this assumes that the supplied boundary conditions are exactly representable as members of the chosen finite-dimensional spaces; more generally, one can always project the boundary conditions into the appropriate finite-dimensional trace spaces. See \eqref{eq: bc NS} and \eqref{eq: bc Maxwell} for the general boundary setting. 

\subsection{Conservation properties} \label{SUBSEC: semi-discrete conservation}
In this section, we analyze some properties of the formulation \eqref{Eq: wf}.
In particular, we show that this semi-discrete formulation preserves conservation of mass and conservation of charge strongly, preserves Gauss' law of magnetism weakly, and preserves the correct energy dissipation rate such that, when no body forces are applied and there is no net flux of energy through the boundary, it conserves the total energy in the ideal limit.

\subsubsection{Conservation of mass} \label{subsec: mass conservation}
For $\bu_{h}\in D(\Omega)$, \eqref{wf c} secures pointwise strong mass conservation, i.e. $\nabla\cdot\bu_{h}=0$, because of the fact that $\nabla\cdot$ maps $D(\Omega)$ surjectively onto $S(\Omega)$, see \eqref{eq: discrete de Rham complex}.

\subsubsection{Conservation of charge} \label{subsub: charge conservation}
Since we have selected $\bH_h\in C(\Omega)$, from \eqref{eq: discrete de Rham complex}, we know that we can always find $\bj_h=\nabla\times\bH_h \in D(\Omega)$ which is pointwise divergence-free; conservation of charge is preserved by the formulation \eqref{Eq: wf}.


\subsubsection{Weak conservation of Gauss’s law for magnetism}\label{subsec: divB conservation}
        \revAPalha{One of the most common forms of the MHD equations is the conservative form (or divergence form); see, for example, \cite{dao_structure_2024}. The discretization of the conservative form poses several challenges. One of these challenges is the involution constraint $\nabla\cdot\boldsymbol{H} = 0$ (or $\nabla\cdot\boldsymbol{B} = 0$), which has been widely discussed in the literature, e.g. \cite{TOTH2000605, Evans1988SimulationOM, wu_provably_2018,brackbill_effect_1980, ramshaw_method_1983, balsara_staggered_1999}. The main root of this challenge is that this formulation is only valid if the involution constraint is satisfied \cite{dao_structure_2024, brackbill_effect_1980}. Moreover, in conservation form, the involution constraint is not an implication of these equations and, therefore, must be enforced. On the other hand, the formulation we employ in this work is self-consistent with regard to the involution constraint. By self-consistent, we mean the following: The involution constraint must only be imposed on the initial condition. Once this condition is satisfied by the initial condition, the evolution equation enforces it at all time instants. As long as numerical discretization preserves a discrete de Rham complex, this equation will automatically be satisfied and the involution constraint will be enforced if the initial condition is divergence-free; see \eqref{eq: divH}.

        In a different context, it has been reported that the formulation leads to better results than the conservation form \cite{brackbill_effect_1980}, when not enforcing the involution constraint. More recently, \cite{dao_structure_2024}, a similar weak formulation has been proposed. In that work, as in this, the weak divergence-free constraint is preserved throughout the whole evolution. As seen in the next section (and in the numerical results), this does not lead to a loss of energy conservation, in contrast to approaches for the conservative form that do not strongly satisfy the divergence-free constraint.}
        
	Since \eqref{wf d} is valid for all $\boldsymbol{b}_{h} \in C_{\boldsymbol{0}}^{\parallel}\left(\Omega,\Gamma_{\widehat{\bH}}\right)$, let us choose $\boldsymbol{b}_h = \nabla q_h$ for an arbitrary $q_h \in G_0(\Omega, \Gamma_{\widehat{\bH}})$. Then, we have
	\begin{equation}
		\ip{\partial_{t}\bH_{h}}{\nabla q_{h}}{\Omega} + \Rm^{-1}\ip{\nabla\times \bH_{h}}{\nabla\times\nabla q_{h}}{\Omega} -\ a\left(\bu_{h}, \bH_{h}, \nabla\times\nabla q_{h}\right) = \ip{\widehat{\bE}}{\mathcal{T}_{\parallel}\nabla q_{h}}{\Gamma_{\widehat{\bE}}}\;.
	\end{equation}
	The second and third terms on the left-hand side are zero because we strongly satisfy $\nabla\times\nabla q_{h} = 0$ since our discrete function spaces constitute a finite-dimensional de Rham complex. We are left with the following equation,
	\begin{equation}
		\ip{\partial_{t}\bH_{h}}{\nabla q_{h}}{\Omega} = \ip{\widehat{\bE}}{\mathcal{T}_{\parallel}\nabla q_{h}}{\Gamma_{\widehat{\bE}}}\,. \label{eq:ampere_law_nabla_q}
	\end{equation}
	Using the fact that $\nabla q_h \in C_{\boldsymbol{0}}^{\parallel}\left(\Omega,\Gamma_{\widehat{\bH}}\right)$, we have
	\begin{equation}
 \begin{aligned}
		\ip{\partial_{t}\bH_{h}}{\nabla q_{h}}{\Omega} = \left\langle \widehat{\boldsymbol{E}}, \mathcal{T}_\parallel \nabla q_{h}\right\rangle_{\Gamma_{\widehat{\boldsymbol{E}}}} &= \left\langle \widehat{\boldsymbol{E}}, \mathcal{T}_\parallel \nabla q_{h}\right\rangle_{\partial\Omega} \\&= -\left\langle \left(\nabla\cdot\right)_{\partial\Omega} \widehat{\boldsymbol{E}}, \mathcal{T}q_{h}\right\rangle_{\partial\Omega} = -\left\langle \left(\nabla\cdot\right)_{\Gamma_{\widehat{\boldsymbol{E}}}} \widehat{\boldsymbol{E}}, \mathcal{T}q_{h}\right\rangle_{\Gamma_{\widehat{\boldsymbol{E}}}}\;.
  \end{aligned}
	\end{equation}
	Finally, using the relation between $\bH \cdot \bn$ and $\widehat{\boldsymbol{E}}$ from \eqref{eq:H_perp_E_tan}, we can rewrite the rightmost term to obtain the following equation,
	\begin{equation}
		\dfrac{\mathrm{d}}{\mathrm{d}t} \left[
			\ip{\bH_{h}}{\nabla q_{h}}{\Omega} - \left\langle \bH \cdot \bn, \mathcal{T}q_{h}\right\rangle_{\Gamma_{\widehat{\boldsymbol{E}}}}
		\right] = 0\;.
	\end{equation}
	The expression inside the square brackets is precisely the weak divergence-free constraint on $\bH_h$. The above equation thus states that the time-derivative of this constraint on $\bH_h$ is zero, and implies that $\bH_{h}$ will remain weakly divergence-free for all $t$ if the initial condition $\bH^0_{h}$ is weakly divergence-free.

\subsubsection{Energy conservation} \label{subsec: energy}
The (semi-)discrete (total) energy is defined as
\begin{equation}\label{eq: total energy}
	\mathcal{E}_{h} := \mathcal{K}_{h}+ \mathcal{M}_{h}
\end{equation}
whose $\mathcal{K}_{h}:= \dfrac{1}{2}\ip{\bu_{h}}{\bu_{h}}{\Omega}$ is the kinetic energy of the fluid and $\mathcal{M}_{h}:=\dfrac{\mathsf{c}}{2}\ip{\bH_{h}}{\bH_{h}}{\Omega}$ is the magnetic energy. Suppose
\begin{equation}
    \hat{u} = 0\;,\quad
    \Gamma_{\widehat{\bu}}=\Gamma_{\widehat{\bE}}=\partial\Omega\;,\quad
    \Gamma_{\widehat{\bw}} = \Gamma_{\widehat{\bH}} = \emptyset\;.
\end{equation}
By replacing $\bv_{h}$ in \eqref{wf a} with $\bu_{h}$, we can obtain
\begin{equation} \label{eq: K energy 1}
	\partial_{t}\mathcal{K}_{h}
+ \Rn^{-1}\ip{\nabla\times\bw_{h}}{\bu_{h}}{\Omega} 
- \mathsf{c}\ 
a\left(\boldsymbol{j}_{h},\bH_{h},\bu_{h}\right)
  = \ip{\bF}{\bu_{h}}{\Omega} - \ip{\widehat{P}}{\mathcal{T}\bu_{h}}{\Gamma_{\widehat{P}}},
\end{equation}
where $\boldsymbol{j}_{h}=\nabla\times\bH_{h} \in D(\Omega)$, and the fluid advection term $a\left(\bw_{h}, \bu_{h}, \bu_{h}\right)$ and the total pressure term $ \ip{P_{h}}{\nabla\cdot\bu_{h}}{\Omega}$ have vanished because of the skew-symmetry of the trilinear \eqref{eq: trilinear} and $\nabla\cdot\bu_{h} = 0$, respectively. Meanwhile, because \eqref{wf b} holds $\forall \boldsymbol{w}_{h}\in C(\Omega)$ and $\bw_{h}\in C(\Omega)$, we know that 
\[-\ip{\bu_{h}}{\nabla\times\bw_{h}}{\Omega}+\ip{\bw_{h}}{\bw_{h}}{\Omega} = - \ip{\widehat{\bu}}{\mathcal{T}_{\parallel}\bw_{h}}{\partial\Omega}.\]
As a results, \eqref{eq: K energy 1} can be expressed as
\begin{equation} \label{eq: K energy 2}
	\partial_{t}\mathcal{K}_{h}
	+ \Rn^{-1}\ip{\bw_{h}}{\bw_{h}}{\Omega} 
	- \mathsf{c}\ 
	a\left(\boldsymbol{j}_{h},\bH_{h},\bu_{h}\right)
	= \ip{\bF}{\bu_{h}}{\Omega} - \ip{\widehat{P}}{\mathcal{T}\bu_{h}}{\Gamma_{\widehat{P}}}
	+ \ip{\widehat{\bu}}{\mathcal{T}_{\parallel}\bw_{h}}{\partial\Omega}.
\end{equation}
Similarly, replacing $\boldsymbol{b}_{h}$ in \eqref{wf d} with $\bH_{h}$ leads to
\begin{equation} \label{eq: B energy 1}
	\dfrac{1}{\mathsf{c}}\partial_{t}\mathcal{M}_{h} + \Rm^{-1}\ip{\boldsymbol{j}_{h}}{\boldsymbol{j}_{h}}{\Omega} 
-a\left(\bu_{h}, \bH_{h}, \boldsymbol{j}_{h}\right)
= \ip{\widehat{\bE}}{\mathcal{T}_{\parallel}\bH _{h}}{\partial\Omega}.
\end{equation}
With \eqref{eq: K energy 2} and \eqref{eq: B energy 1}, we can finally derive the rate of change of the discrete energy,
\begin{equation} \label{eq: time derivative of energy}
	\begin{aligned}
		\partial_{t} \mathcal{E}_{h} 
		&= \partial_{t}\mathcal{K}_{h} + \partial_{t}\mathcal{M}_{h}\\
		&= \underbrace{\ip{\bF}{\bu_{h}}{\Omega}}_{\text{(i)}}\ 
		\underbrace{- \Rn^{-1}\mathcal{S}_{h}}_{\text{(ii)}}\ 
		\underbrace{- \mathsf{c}\Rm^{-1}\mathcal{J}_{h}}_{\text{(iii)}}
		\underbrace{- \ip{\widehat{P}}{\mathcal{T}\bu_{h}}{\Gamma_{\widehat{P}}}}_{\text{(iv)}}
		+ \underbrace{\ip{\widehat{\bu}}{\mathcal{T}_{\parallel}\bw_{h}}{\partial\Omega} }_{\text{(v)}}
		+\underbrace{\mathsf{c}\ip{\widehat{\bE}}{\mathcal{T}_{\parallel}\bH _{h}}{\partial\Omega}} _{\text{(vi)}},
	\end{aligned}
\end{equation}
where 
$\mathcal{S}_{h} := \ip{\bw_{h}}{\bw_{h}}{\Omega} $,
$\mathcal{J}_{h} := \ip{\boldsymbol{j}_{h}}{\boldsymbol{j}_{h}}{\Omega}$. The two nonlinear terms have canceled each other due to the skew-symmetry of the trilinear form \eqref{eq: trilinear}; the exchange of energy between $\mathcal{K}_{h}$ and $\mathcal{M}_{h}$ is exactly captured.

In \eqref{eq: time derivative of energy}, (i) represents the amount of mechanical work done on the system per unit time, (ii) and (iii) are dissipation rates due to viscosity and electric resistance, respectively, and terms (iv) - (vi) are the net flux of energy through the domain boundary.
Clearly, if $\bF=\boldsymbol{0}$, when the MHD flow is ideal (i.e. $\Rn=\Rm=\infty$) and there is no net flux of energy through the domain boundary, formulation \eqref{Eq: wf} preserves energy, i.e.
\[
\partial_{t} \mathcal{E}_{h}  =  \partial_{t}\mathcal{K}_{h} + \partial_{t}\mathcal{M}_{h} = 0.
\]

\section{Fully discrete weak formulation} \label{Sec: discretization}

If a conservative temporal discrete scheme, for example, the Crank-Nicolson scheme, is applied to the semi-discrete formulation \eqref{Eq: wf} directly, we can obtain a fully discrete formulation of the same conservation and dissipation properties (except that the time derivative of energy is now expressed in a temporally discrete format). See \ref{App: coupled discretization}.
The proof is straightforward; it mimics the analyses in Section~\ref{SUBSEC: semi-discrete conservation}.
Therefore, the rest of this paper is devoted to a decoupled formulation with computational efficiency in mind.

\subsection{Decoupled temporal discretization}
Consider a set of time instants
\begin{equation}\label{eq: time sequence}
\left\lbrace t^{0}, t^{\frac{1}{2}}, t^{1}, t^{1+\frac{1}{2}}, t^{2}, \cdots \right\rbrace,
\end{equation}
where $t^{0} = 0 $, $\Delta t = t^{k}-t^{k-1} > 0$, $k\in\left\lbrace1,2,3,\cdots\right\rbrace$, is a constant time interval, and $t^{k-\frac{1}{2}} = \dfrac{t^{k-1} + t^{k}}{2}$. We use a superscript to denote the evaluation of a variable at a particular time instant. For example, $\bu^{k}_{h}:= \bu_{h}\left(\boldsymbol{x}, t^{k}\right)$. 

A fully discrete, temporally decoupled formulation of \eqref{Eq: wf} (hereinafter referred to as the decoupled formulation) is expressed as follows. Given $\bF\in\left[L^2(\Omega)\right]^3$, natural boundary conditions 
$\widehat{P}\in H^{1/2}\left(\Omega;\Gamma_{\widehat{P}}\right),\ \widehat{\bu}\in\mathcal{T}H\left(\mathrm{curl};\Omega,\Gamma_{\widehat{\bu}}\right),\  \widehat{\bE}\in \mathcal{T}H\left(\mathrm{curl};\Omega,\Gamma_{\widehat{\boldsymbol{E}}}\right)$,
and initial conditions $\left(\bu_{h}^{0}, \bw_{h}^{0}, \bH^{\frac{1}{2}}_{h}\right)\in D\left(\Omega\right)\times C\left(\Omega\right)\times C\left(\Omega\right)$, for time-step index $k=1,2,3,\cdots $ successively,
\vspace{0.2cm}\\\noindent {(\textbf{step 1})} seek $\left(\boldsymbol{u}^{k}_{h}, \bw^{k}_{h}, P^{k-\frac{1}{2}}_{h} \right)\in D_{\hat{u}}(\Omega,\Gamma_{\hat{u}})\times C_{\widehat{\bw}}^{\parallel}(\Omega,\Gamma_{\widehat{\bw}})\times S(\Omega)$, such that,
$\forall \left(\bv_{h},\boldsymbol{w}_{h}, q_{h}\right)\in D_{0}\left(\Omega,\Gamma_{\hat{u}}\right) \times C^{\parallel}_{\boldsymbol{0}}\left(\Omega,\Gamma_{\widehat{\bw}}\right)\times S\left(\Omega\right)$, 
\begin{subequations}\label{Eq: wf step 1}
	\begin{align}
		\label{wf step1 a}
		\left\langle\dfrac{\bu^{k}_{h}-\bu^{k-1}_{h}}{\Delta t}, \bv_{h}\right\rangle_{\Omega}
		+
		a\left(\dfrac{\bw^{k-1}_{h} + \bw^{k}_{h}}{2}, \dfrac{\bu^{k-1}_{h} + \bu^{k}_{h}}{2}, \bv_{h}\right)
		&+ \Rn^{-1}\ip{\nabla\times\dfrac{\bw^{k-1}_{h} + \bw^{k}_{h}}{2}}{\bv_{h}}{\Omega} \\ 
		\nonumber
		- \mathsf{c}\ 
		a\left(\nabla\times\bH^{k-\frac{1}{2}}_{h},\bH^{k-\frac{1}{2}}_{h},\bv_{h}\right)
		- \ip{P^{k-\frac{1}{2}}_{h}}{\nabla\cdot\bv_{h}}{\Omega}  &= \ip{\bF^{k-\frac{1}{2}}}{\bv_{h}}{\Omega} - \ip{\widehat{P}^{k-\frac{1}{2}}}{\mathcal{T}\bv_{h}}{\Gamma_{\widehat{P}}},\\
		\label{wf step1 b}
		-\ip{\bu^{k}_{h}}{\nabla\times\boldsymbol{w}_{h}}{\Omega}+\ip{\bw^{k}_{h}}{\boldsymbol{w}_{h}}{\Omega} & = - \ip{\widehat{\bu}^{k}}{\mathcal{T}_{\parallel}\boldsymbol{w}_{h}}{\Gamma_{\widehat{\bu}}} ,\\
		\label{wf step1 c}
		\ip{\nabla\cdot\bu^{k}_{h}}{q_{h}}{\Omega} & = 0,
	\end{align}
\end{subequations}
\vspace{0.2cm}\\\noindent{(\textbf{step 2})} seek $\bH^{k+\frac{1}{2}}_{h}\in C^{\parallel}_{\widehat{\bH}}(\Omega,\Gamma_{\widehat{\bH}}) $, such that, 
$\forall \boldsymbol{b}_{h}\in  C_{\boldsymbol{0}}^{\parallel}\left(\Omega,\Gamma_{\widehat{\bH}}\right)$,
\begin{equation}\label{Eq: wf step 2}
		\begin{aligned}
			\ip{\dfrac{\bH^{k+\frac{1}{2}}_{h} - \bH^{k-\frac{1}{2}}_{h}}{\Delta t}}{\boldsymbol{b}_{h}}{\Omega} + \Rm^{-1}\ip{\nabla\times \dfrac{\bH^{k-\frac{1}{2}}_{h} + \bH^{k+\frac{1}{2}}_{h}}{2}}{\nabla\times\boldsymbol{b}_{h}}{\Omega} 
			\qquad &\\
			-\ a\left(\bu^{k}_{h}, \dfrac{\bH^{k-\frac{1}{2}}_{h} + \bH^{k+\frac{1}{2}}_{h}}{2}, \nabla\times\boldsymbol{b}_{h}\right)
			&= \ip{\widehat{\bE}^{k}}{\mathcal{T}_{\parallel}\boldsymbol{b}_{h}}{\Gamma_{\widehat{\bE}}}.
		\end{aligned}
\end{equation}

\begin{figure}[h!]
	\centering
	\includegraphics[width=0.7\linewidth]{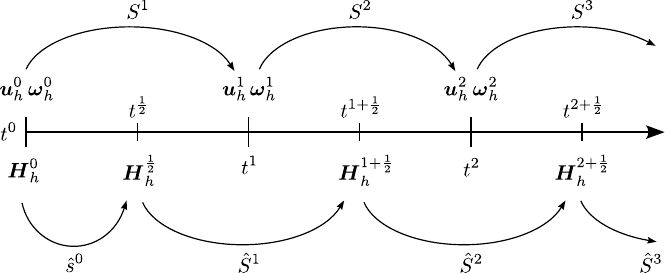}
	\caption{An illustration of the leapfrog temporal scheme for the decoupled formulation. The time-step sequence is $\hat{s}^{0}\to S^{1}\to\hat{S}^1\to S^2\to \hat{S}^2\to\cdots$. Steps $S^{1}, S^{2},S^3, \cdots$ refer to the step 1, see \eqref{Eq: wf step 1}, of the decoupled formulation, and steps $ \hat{S}^{1}, \hat{S}^{2},\hat{S}^3, \cdots$ refer to the step 2, see \eqref{Eq: wf step 2}, of the decoupled formulation. The pre-step $\hat{s}^{0}$ computes $\bH_{h}^{\frac{1}{2}}$. A similar scheme has been used in a dual-field discretization of incompressible Navier-Stokes equations, see \cite[Fig.~1]{ZHANG2022110868}.}
	\label{fig:scheme illustration}
\end{figure}

Note that the decoupled formulation takes $\bH_{h}^{\frac{1}{2}}$, instead of $\bH_{h}^{0}$, as an initial condition. Thus, to initiate the iterations, we need to compute $\bH_{h}^{\frac{1}{2}}$ by applying, for example, an explicit Euler scheme or the Crank-Nicolson scheme as discussed in Section~\ref{Sec: discretization} to \eqref{Eq: wf} at the half-time-step from $t^{0}$ to $t^{\frac{1}{2}}$. 
An illustration of the overall temporal scheme for the decoupled formulation is shown in Fig.~\ref{fig:scheme illustration}.

This leapfrog-type temporal scheme is analogous to the symplectic Störmer–Verlet leapfrog integrator \cite{Hairer_Lubich_Wanner_2003} in Hamiltonian mechanics. Seeing that the decoupled formulation has linearized the Maxwell part, see \eqref{Eq: wf step 2}, it does not fully linearize the Navier-Stokes part. More specifically, the second term of \eqref{wf step1 a}, namely, the fluid convection term, is yet to be left non-linear. Thus, a nonlinear system solver, for example, the Picard method or the Newton-Raphson method, still needs to be employed.

\subsection{Conservation properties}\label{Subsec: properties temporally discrete}
Clearly, the formulation \eqref{Eq: wf step 1} strongly preserves conservation of mass, i.e., $\nabla\cdot\bu_{h}^{k}= 0$, and conservation of charge, i.e., $\nabla\cdot\bj_h^{k+\frac{1}{2}}=0$, pointwise, and weakly preserves Gauss’s law for magnetism, for same reasons as in Section~\ref{subsec: mass conservation}, Section~\ref{subsub: charge conservation} and Section~\ref{subsec: divB conservation}, respectively; the temporal discretization scheme maintains these conservation properties at the fully discrete level.

For energy conservation, it has been shown that, for the semi-discrete formulation, a key to deriving it is the canceling of the two nonlinear terms, see \eqref{eq: K energy 2} - \eqref{eq: time derivative of energy}. We now study whether a similar canceling happens for the decoupled formulation. And, to simplify the notations, we use a periodic domain, i.e., whose $\partial\Omega=\emptyset$ or vanishing boundary conditions, i.e. $\widehat{P}= 0$ and $\widehat{\bu}=\widehat{\boldsymbol{E}}=\boldsymbol{0}$ on $\partial\Omega\times(0,T]$. Throughout the paper, we use the midpoint rule, i.e., for example, 
\[\bu_{h}^{k-\frac{1}{2}}:=\dfrac{\bu_{h}^{k-1}+\bu_{h}^{k}}{2}.\]

If we replace $\bv_{h}$ in \eqref{wf step1 a} with $\bu_{h}^{k-\frac{1}{2}}$, by performing similar analyses as those in Section~\ref{subsec: energy}, we can find that
\begin{equation*} \label{eq: K energy fd}
	\begin{aligned}
		\left\langle\frac{\bu^{k}_{h}-\bu^{k-1}_{h}}{\Delta t},\frac{\bu^{k-1}_{h}+\bu^{k}_{h}}{2}\right\rangle_{\Omega}
		+ \Rn^{-1}\ip{\bw_{h}^{k-\frac{1}{2}}}{\bw_{h}^{k-\frac{1}{2}}}{\Omega}
		- \mathsf{c}\ a\left(\boldsymbol{j}^{k-\frac{1}{2}}_{h},\bH^{k-\frac{1}{2}}_{h},\bu^{k-\frac{1}{2}}_{h}\right)
		= \ip{\bF^{k-\frac{1}{2}}}{\bu^{k-\frac{1}{2}}_{h}}{\Omega},
	\end{aligned}
\end{equation*}
where $\bj_{h}^{k-\frac{1}{2}} = \nabla\times \bH_{h}^{k-\frac{1}{2}}$. This leads to, if $\mathcal{K}^{k} := \dfrac{1}{2}\ip{\bu^{k}}{\bu^{k}_{h}}{\Omega} $ denotes the discrete kinetic energy,
\begin{equation}\label{eq: K}
	\begin{aligned}
		\frac{ \mathcal{K}_{h}^{k}- \mathcal{K}_{h}^{k-1}}{\Delta t}
		= \mathcal{F}^{k-\frac{1}{2}}
		- \Rn^{-1}\mathcal{S}_{h}^{k-\frac{1}{2}}
		+ \mathsf{c}\mathcal{A}_{h}^{k-\frac{1}{2}},
	\end{aligned}
\end{equation}
where 
$\mathcal{F}^{k-\frac{1}{2}} := \ip{\bF^{k-\frac{1}{2}}}{\bu^{k-\frac{1}{2}}_{h}}{\Omega} $,
$\mathcal{S}_{h}^{k-\frac{1}{2}} := \ip{\bw_{h}^{k-\frac{1}{2}}}{\bw_{h}^{k-\frac{1}{2}}}{\Omega} $ and
$\mathcal{A}_{h}^{k-\frac{1}{2}}:=a\left(\boldsymbol{j}^{k-\frac{1}{2}}_{h},\bH^{k-\frac{1}{2}}_{h},\bu^{k-\frac{1}{2}}_{h}\right)$. 
Similarly, by replacing $\boldsymbol{b}_{h}$ in \eqref{Eq: wf step 2} with $\bH^{k}_{h}$, we can obtain
\begin{equation}\label{eq: step k}
	\frac{\mathcal{M}_{h}^{k+\frac{1}{2}}- \mathcal{M}_{h}^{k-\frac{1}{2}}}{\Delta t}= -\mathsf{c}\Rm^{-1}\mathcal{J}_{h}^{k} - \mathsf{c}\mathcal{A}^{k}_{h},
\end{equation}
where $\mathcal{M}_{h}^{k+\frac{1}{2}}=\dfrac{\mathsf{c}}{2}\ip{\bH^{k+\frac{1}{2}}}{\bH^{k+\frac{1}{2}}}{\Omega}$ denotes the discrete magnetic energy, and $\mathcal{J}_{h}^{k}:=\ip{\boldsymbol{j}^{k}_{h}}{\boldsymbol{j}^{k}_{h}}{\Omega} $. Repeating this analysis for the time-step indexed $(k-1)$ gives
\begin{equation}\label{eq: step k - 1}
	\frac{ \mathcal{M}_{h}^{k-\frac{1}{2}}- \mathcal{M}_{h}^{k-\frac{3}{2}}}{\Delta t} = -\mathsf{c}\Rm^{-1}\mathcal{J}_{h}^{k-1} - \mathsf{c}\mathcal{A}^{k-1}_{h}.
\end{equation}
\eqref{eq: step k} and \eqref{eq: step k - 1} together shows
\begin{equation}\label{eq: B}
	\frac{ \tilde{\mathcal{M}}_{h}^{k}- \tilde{\mathcal{M}}_{h}^{k-1}}{\Delta t} = 
-\mathsf{c}\Rm^{-1}\tilde{\mathcal{J}}_{h}^{k-\frac{1}{2}} 
- \mathsf{c}\tilde{\mathcal{A}}^{k-\frac{1}{2}}_{h},
\end{equation}
where we have used the tilde sign to denote that it is an average of two integral quantities, i.e., for example, 
\begin{equation} \label{eq p e}
\tilde{\mathcal{M}}_{h}^{k} := \dfrac{\mathcal{M}_{h}^{k-\frac{1}{2}} + \mathcal{M}_{h}^{k+\frac{1}{2}}}{2}
\end{equation}
which differs from ${\mathcal{M}}_{h}^{k} =  \dfrac{\mathsf{c}}{2}\ip{\bH^{k}}{\bH^{k}}{\Omega}.$ 
With \eqref{eq: K} and \eqref{eq: B}, we can conclude that the decoupled formulation dissipates a discrete energy defined as 
\begin{equation}\label{eq: discrete energy definitaion}
	\tilde{\mathcal{E}}_{h}^{k}:= {\mathcal{K}}_{h}^{k} + \tilde{\mathcal{M}}_{h}^{k}
\end{equation}
with a discrete rate
\begin{equation}\label{Eq: discrete energy dissipation}
	\frac{\tilde{\mathcal{E}}_{h}^{k} - \tilde{\mathcal{E}}_{h}^{k-1}}{\Delta t} = \mathcal{F}^{k-\frac{1}{2}}
- \Rn^{-1}\mathcal{S}_{h}^{k-\frac{1}{2}}
- \mathsf{c}\Rm^{-1}\tilde{\mathcal{J}}_{h}^{k-\frac{1}{2}}
+ \mathsf{c}\left( \mathcal{A}_{h}^{k-\frac{1}{2}} -\tilde{\mathcal{A}}^{k-\frac{1}{2}}_{h}\right).
\end{equation}
\revABrugnoli{It is not guarantee that the formulation will dissipate as the last term is undefined;}
when $\boldsymbol{f}=\boldsymbol{0}$ and the flow is ideal, we have
\begin{equation*}\label{Eq: discrete energy conservation}
	\frac{\tilde{\mathcal{E}}_{h}^{k} - \tilde{\mathcal{E}}_{h}^{k-1}}{\Delta t} = \mathsf{c}\left( \mathcal{A}_{h}^{k-\frac{1}{2}} -\tilde{\mathcal{A}}^{k-\frac{1}{2}}_{h}\right).
\end{equation*}
As generally $ \mathcal{A}_{h}^{k-\frac{1}{2}} -\tilde{\mathcal{A}}^{k-\frac{1}{2} }\neq 0$, the decoupled formulation does not conserve the discrete energy \eqref{eq: discrete energy definitaion}. 
Exact preservation of the discrete energy would require a fully implicit scheme coupling the fluid and Maxwell parts. This can be achieved, for instance, using the midpoint rule as the Hamiltonian is quadratic in the velocity and magnetic field \revABrugnoli{and the midpoint rule preserves quadratic invariants \cite{hairer2006geometric}. Many integration strategies have been proposed in recent years to ensure energy preservation without the need to solve nonlinear systems. These schemes rely on the idea of energy quadratization via the introduction of a scalar auxiliary variable and lead to linearly implicit schemes \cite{shen2018sav}. The scalar auxiliary variable approach combines the implicit midpoint scheme with a St\"ormer-Verlet integrator to achieve second-order accuracy in time and energy conservation. For the MHD equations, the nonlinear coupling between fluid and electromagnetic systems makes it difficult to achieve second-order accuracy and energy conservation without the need for a fully implicit nonlinear solver.}

In short, the decoupled formulation breaks the nonlinear MHD formulation into two evolution equations representing the Navier-Stokes part and the Maxwell part, respectively, partially linearizes the former and fully linearizes the latter. It keeps strong conservation of mass, strong conservation of charge, and weak conservation of Gauss’s law for magnetism while sacrificing conservation of energy. \revABrugnoli{More importantly, it allows to keep second-order accuracy in time while confining the nonlinear iterations to the fluid system.} 


\section{Numerical tests}\label{Sec: numerical tests}
In this section, we present some results of numerical tests using the proposed method. The test cases are manufactured solution tests, the Orszag-Tang vortex flow, and magnetic lid-driven cavity flow. We use the mimetic spectral elements \cite{Kreeft2011,zhang2022phd} for all tests. The degree of the basis functions is denoted by $N$. Any other set of finite elements that also satisfies the discrete de Rham Hilbert complex \eqref{eq: discrete de Rham complex} and the $L^2$-regularity for the triple products, see Section~\ref{subseq:finite_dim_setting} works for the proposed method. A combination of Lagrange (continuous Galerkin) elements of degree $N$,  the first kind N{\'{e}}d{\'{e}}lec $H(\mathrm{curl})$-conforming elements of degree $N$ \cite{nedelec},  Raviart–Thomas elements of degree $N$ \cite{raviart606mixed},  and discontinuous Galerkin elements of degree $(N-1)$ \cite{carstensen2016breaking} is a classic candidate \cite{boffi2013mixed,doi:10.1137/1.9781611977738}.

\subsection{Manufactured solution tests}\label{Subsec: manu tests}
Two sets of manufactured solutions are employed for convergence tests and conservation and dissipation tests, respectively.
\subsubsection{Convergence tests} \label{Subsec: convergence}
Suppose three-dimensional manufactured solutions
\[
\bu = \begin{bmatrix}
	\cos(x) \sin(y) \sin(z) e ^ t & \sin(x) \cos(y) \sin(z) e ^ {t} & -2\sin(x) \sin(y) \cos(z) e ^ {t}
\end{bmatrix}^{\mathsf{T}},
\]
\[
\bE = \begin{bmatrix}
	\cos(x) \sin(y) \sin(z/2) e ^ t & \sin(x/2) \cos(y) \sin(z) e ^ t & -\sin(x) \sin(y) \cos(z) e ^ t
\end{bmatrix}^{\mathsf{T}},
\]
\[
P = \cos(x) \cos(y) \cos(z) e ^ {-t}\quad \text{and}\quad
\bH^{0}=\bB^{0} = \boldsymbol{0},
\]
solve \eqref{Eq: dimensionless ast 0}. Manufactured solutions for the remaining variables can be obtained from
$\bw = \nabla\times \bu$,
$\bH=\bB = \int_{0}^{t} \left( -\nabla\times\bE\right)\mathrm{d}t $, $\bj=\nabla\times \bH$. Given $\Rn$ and $\mathsf{c}$, we can get the analytical expression of $\bF$ through \eqref{Eq: dimensionless ast 0 a}. And if $\Rm$ is known, we can find an extra source, denoted by $\boldsymbol{e}$, needed to balance the relation
\[
\boldsymbol{e} =\Rm^{-1}\bj  - \left( \bE  +\bu \times\bB \right).
\]
This will lead to an extra term, $\left\langle\boldsymbol{e},\nabla\times\boldsymbol{b}_{h}\right\rangle_{\Omega}$, in the right hand sides of, for example, \eqref{wf d} and \eqref{Eq: wf step 2}.

\begin{figure}[!htb]
	\centering
	\begin{minipage}[c]{1\textwidth}
		\centering{
			\subfloat{
				\begin{minipage}[b]{0.5\textwidth}
					\centering
					\includegraphics[width=0.8\linewidth]{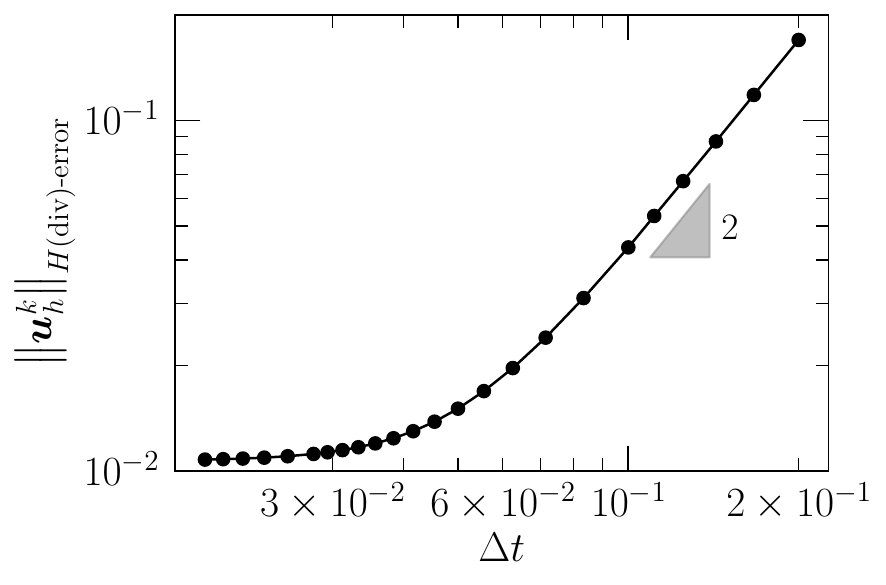}
				\end{minipage}
			}
			\subfloat{
				\begin{minipage}[b]{0.5\textwidth}
					\centering
					\includegraphics[width=0.8\linewidth]{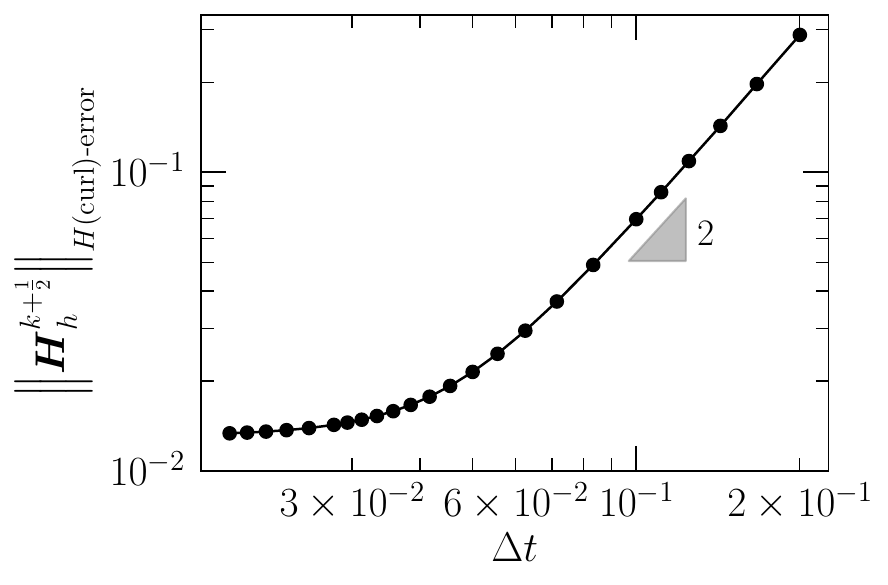}
				\end{minipage}
			}
		}
	\end{minipage}%
	\caption{Results of
		$\left\| \boldsymbol{u}^k_h \right\|_{H(\mathrm{div})\text{-error}}$ and 
		$\left\| \bH ^{k+\frac{1}{2}}_h \right\|_{H(\mathrm{curl})\text{-error}}$
		for the temporal convergence tests at $N=3$, $K=16$, $\Delta t \in \left\lbrace\frac{1}{5},\frac{1}{6},\cdots,\frac{1}{10},\frac{1}{12},\cdots,\frac{1}{36},\frac{1}{40},\cdots,\frac{1}{56}\right\rbrace$.
	}
	\label{fig:temporal convergence tests}
\end{figure}

\begin{figure}[!htb]
	\centering
	\begin{minipage}[c]{1\textwidth}
		\centering{
			\subfloat{
				\begin{minipage}[b]{0.5\textwidth}
					\centering
					\includegraphics[width=0.8\linewidth]{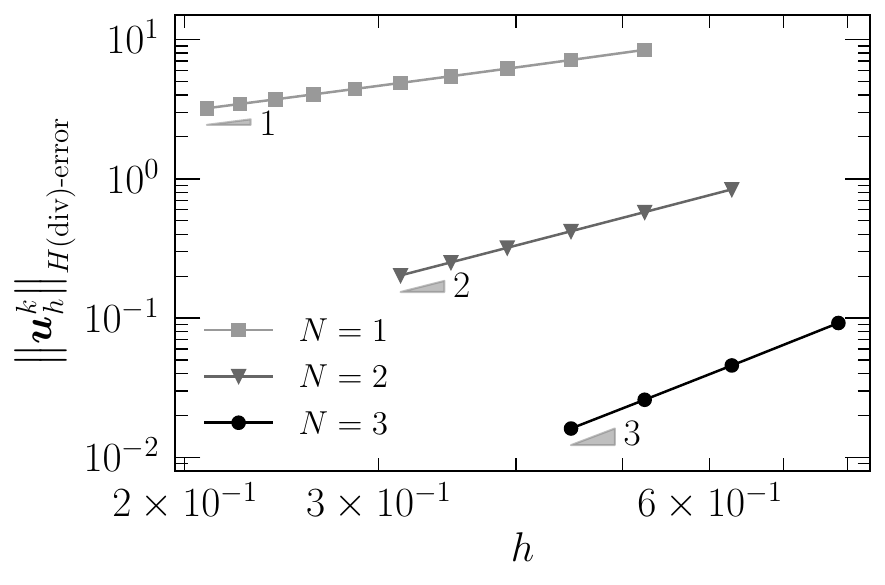}
				\end{minipage}
			}
			\subfloat{
				\begin{minipage}[b]{0.5\textwidth}
					\centering
					\includegraphics[width=0.8\linewidth]{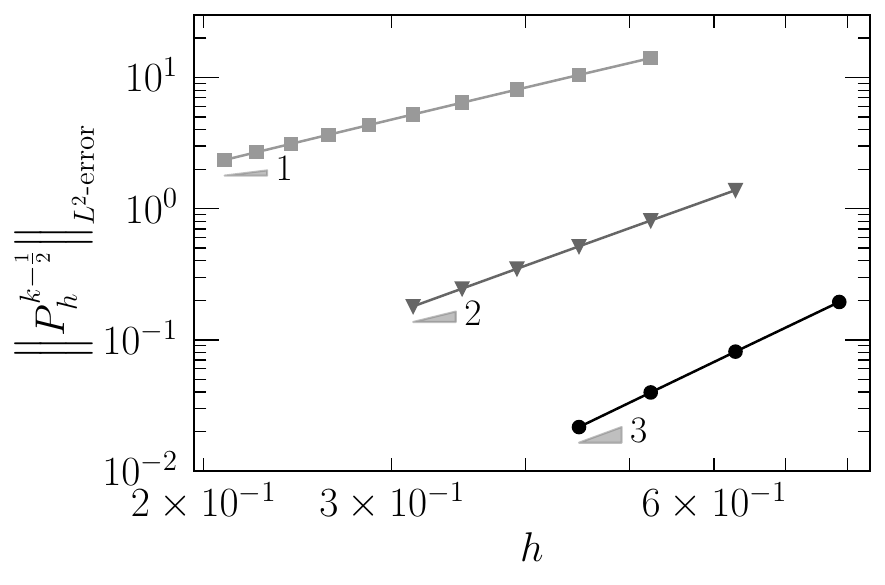}
				\end{minipage}
			}\\
			\subfloat{
				\begin{minipage}[b]{0.5\textwidth}
					\centering
					\includegraphics[width=0.8\linewidth]{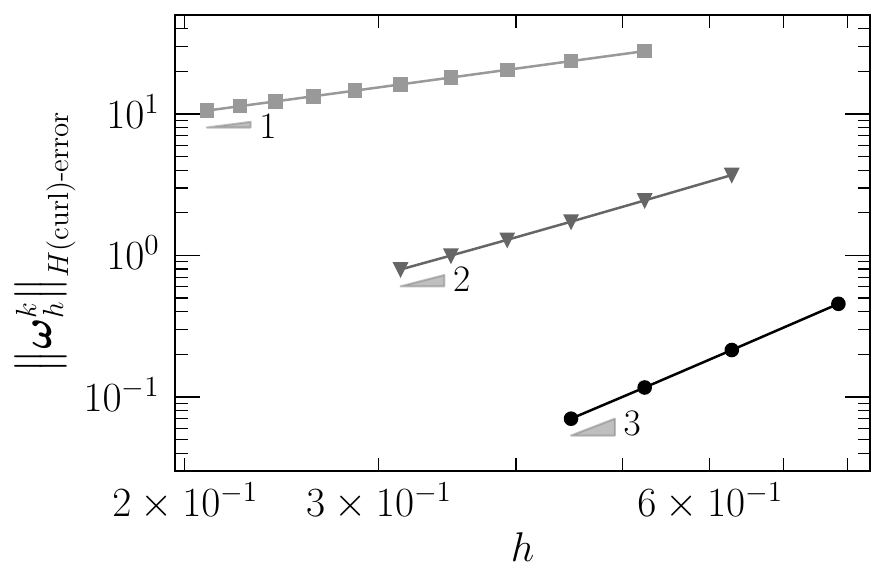}
				\end{minipage}
			}
			\subfloat{
				\begin{minipage}[b]{0.5\textwidth}
					\centering
					\includegraphics[width=0.8\linewidth]{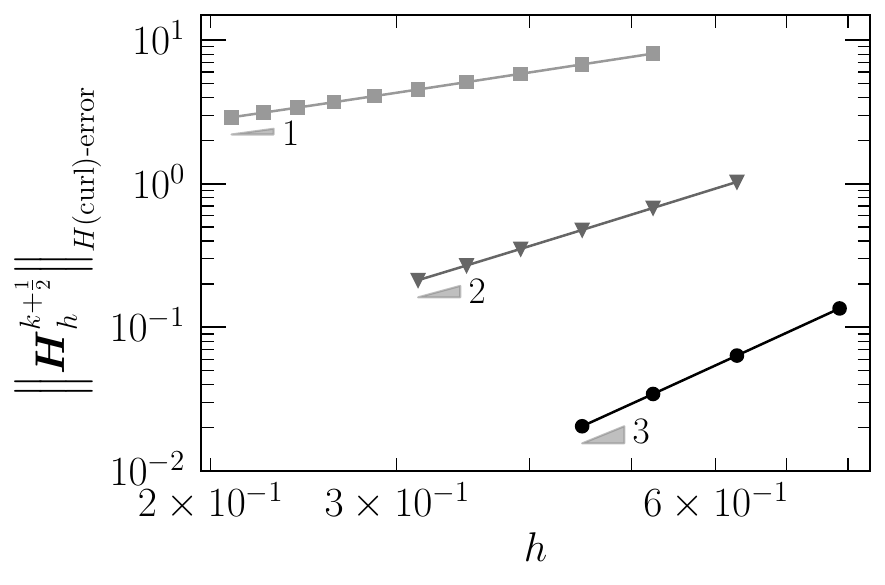}
				\end{minipage}
			}
		}
	\end{minipage}%
	\caption{Results of
		$\left\| \boldsymbol{u}^k_h \right\|_{H(\mathrm{div})\text{-error}}$,
		$\left\| P^{k-\frac{1}{2}}_h\right\|_{L^2\text{-error}}$,
		$\left\| \boldsymbol{\omega}^k_h \right\|_{H(\mathrm{curl})\text{-error}}$ and
		$\left\| \bH ^{k+\frac{1}{2}}_h \right\|_{H(\mathrm{curl})\text{-error}}$
		for the spatial convergence tests. We use $K\in\left\lbrace12,14,\cdots,32\right\rbrace$ for $N=1$, $K\in\left\lbrace10,12,\cdots,20\right\rbrace$ for $N=2$ and $K\in\left\lbrace8,10,\cdots,14\right\rbrace$ for $N=3$. And $\Delta t = \frac{1}{100}$ is employed to avoid pollution of the temporal discretization error, see Fig.~\ref{fig:temporal convergence tests}.
		}
	\label{fig:convergence tests}
\end{figure}

\begin{figure}[h!]
	\centering
	\includegraphics[width=0.45\linewidth]{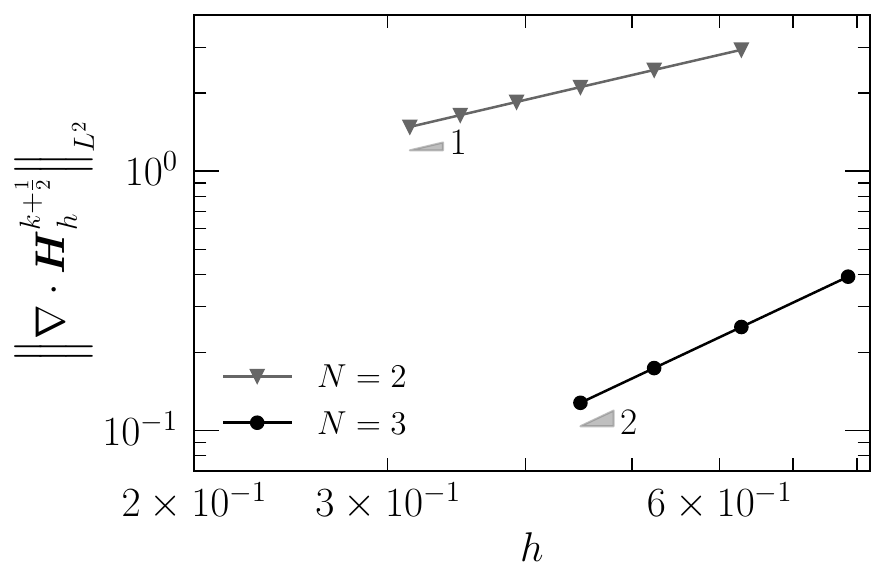}
	\caption{Results of $\left\|\nabla\cdot\boldsymbol{H}^{k+\frac{1}{2}}_h\right\|_{L^2}$ for the spatial convergence tests. We use $K\in\left\lbrace10,12,\cdots,20\right\rbrace$ for $N=2$, $K\in\left\lbrace8,10,\cdots,14\right\rbrace$ for $N=3$, and $\Delta t = \frac{1}{100}$. In this work, the employed finite-dimensional space for $\boldsymbol{H}_h \in C(\Omega)$ is the mimetic spectral space which is a polynomial space $\mathcal{P}^{N-1,N,N}\times\mathcal{P}^{N,N-1,N}\times\mathcal{P}^{N,N,N-1}$ on orthogonal meshes. So we can compute exact divergence of $\boldsymbol{H}_h$ element-wise, and, when $N=1$, $\nabla\cdot\boldsymbol{H}_h\equiv 0$. For the present results, $\left\|\nabla\cdot\boldsymbol{H}^{k+\frac{1}{2}}_h\right\|_{L^2}:=\sqrt{\sum_i H_i^2}$ where $H_i$ is the $L^2$-norm of $\nabla\cdot\boldsymbol{H}_h^{k+\frac{1}{2}}$ in the $i$th element.}
	\label{fig:divH convergence rate}
\end{figure}

For the convergence tests, the parameters are set to $\Rn = 1$, $\Rm=1$, $\mathsf{c}=1$. And the domain is selected to be $\Omega = \left(x, y, z\right)\in \left[0, 2\pi\right]^3$ whose boundary is partitioned, see \eqref{eq: bc NS} and \eqref{eq: bc Maxwell}, such that
\[
\Gamma_{\widehat{P}} = \Gamma_{x}^{-} \cup \Gamma_{y}^{+} \cup \Gamma_{z}^{+},
\]
\[
\Gamma_{\widehat{\bu}} = \Gamma_{x}^{+} \cup \Gamma_{y}^{-} \cup \Gamma_{z}^{+},
\]
\[
\Gamma_{\widehat{\bE}} = \Gamma_{x}^{+} \cup \Gamma_{y}^{+} \cup \Gamma_{z}^{-},
\]
where, for example, $\Gamma_{x}^{-}$ represents the face $\left(x,y,z\right) \in 0\times(0,2\pi)\times(0,2\pi)$ and $\Gamma_{x}^{+}$ represents the face $\left(x,y,z\right) \in 2\pi\times(0,2\pi)\times(0,2\pi)$. Let $K$ and $\kappa$ be two positive integers. We generate a mesh of $K^3$ uniform cubic elements in $\Omega$, and $h=\frac{2\pi}{K}$ is the size, i.e. the edge length, of elements. On this mesh, given initial conditions, 
boundary conditions,  $\boldsymbol{f}$ and $\boldsymbol{e}$ according to the manufactured solutions, the decoupled formulation, \eqref{Eq: wf step 1} and \eqref{Eq: wf step 2}, is solved with a constant time-step interval $\Delta t = \frac{1}{\kappa}$. Errors between simulation and manufactured solutions are measured at $t=1$, i.e. when time-step index $k=\kappa$.
Results in Fig.~\ref{fig:temporal convergence tests} show that the decoupled temporal discretization has a second-order accuracy for both evolution equations. In Fig.~\ref{fig:convergence tests}, optimal spatial convergence rates are observed for all variables. Results supporting that the formulation weakly preserves Gauss' law for magnetism are presented in Fig.~\ref{fig:divH convergence rate}.

 \subsubsection{Conservation and dissipation tests} \label{Subsec: energy-dissipation}
\begin{figure}[h!]
	\centering
	\subfloat{
		\begin{minipage}[b]{0.5\textwidth}
			\centering
			\includegraphics[width=0.88\linewidth]{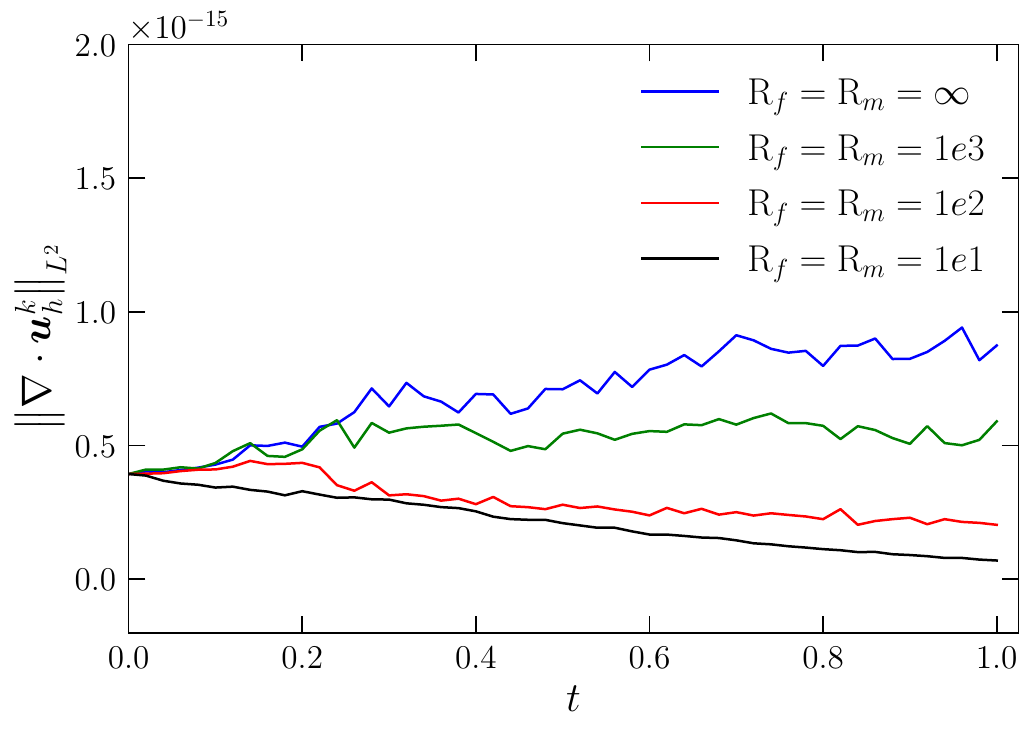}
		\end{minipage}
	}
\\
	
	\begin{minipage}[c]{0.5\textwidth}
		\centering{
			\subfloat{
				\begin{minipage}[b]{1\textwidth}
					\centering
					\includegraphics[width=1\linewidth]{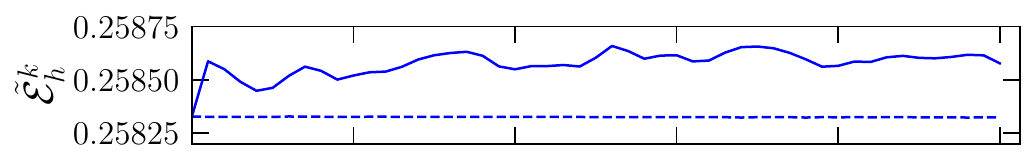}
				\end{minipage}
			}\\\vspace{-0.5cm}
			\subfloat{
				\begin{minipage}[b]{1\textwidth}
					\centering
					\includegraphics[width=1\linewidth]{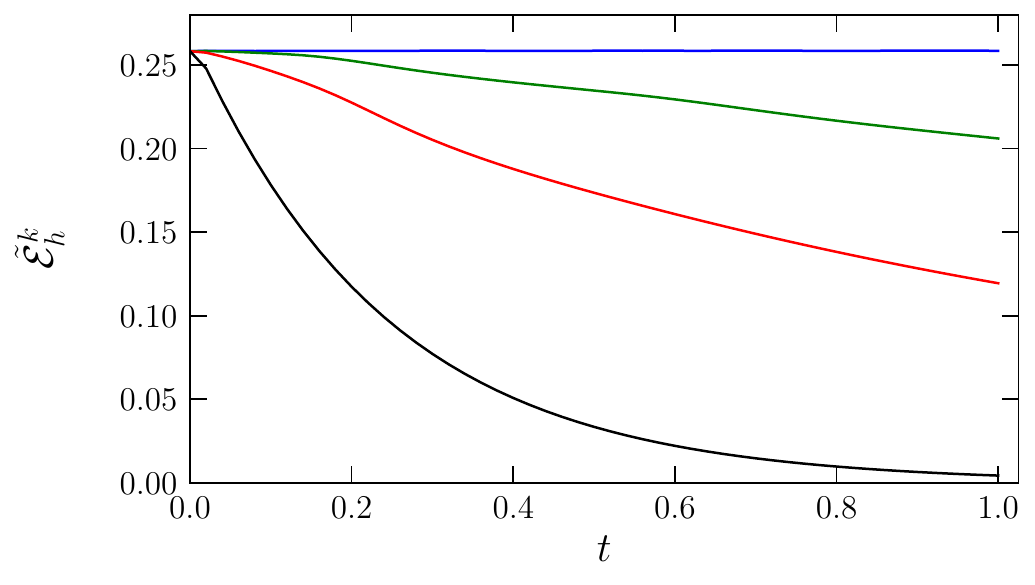}
				\end{minipage}
			}
		}
	\end{minipage}%
	\begin{minipage}[c]{0.5\textwidth}
	\centering{
		\subfloat{
			\begin{minipage}[b]{0.9\textwidth}
				\centering
				\includegraphics[width=1\linewidth]{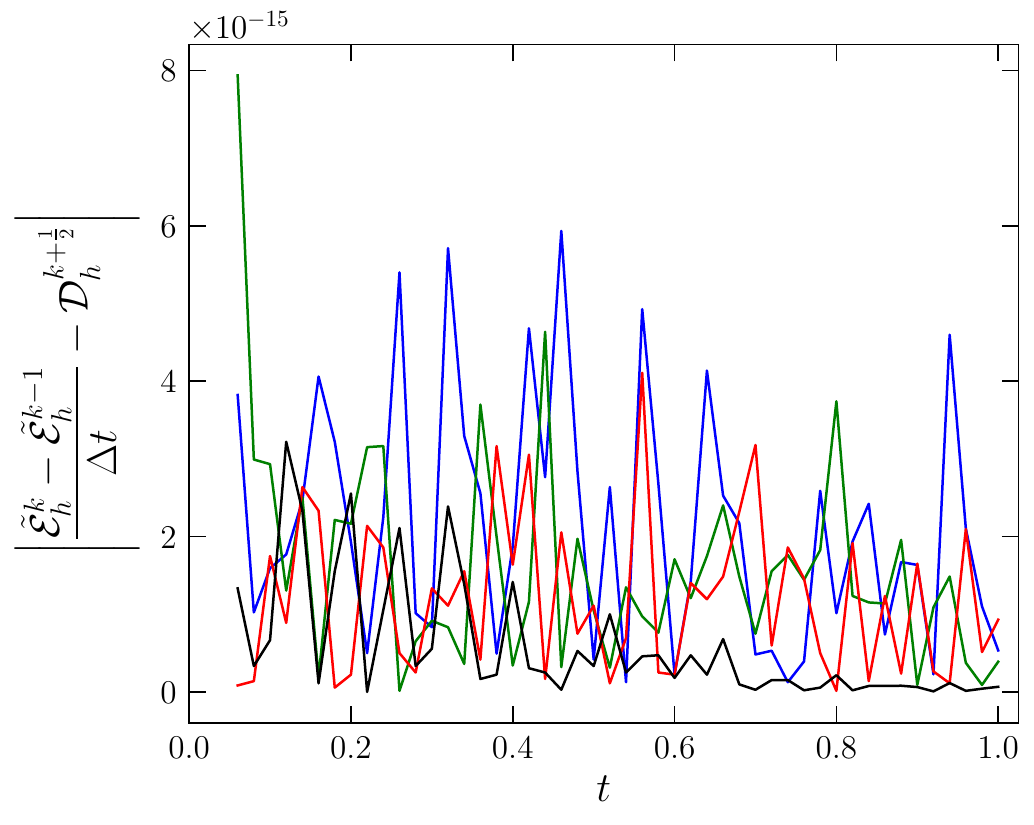}
			\end{minipage}
		}\\\vspace{0.15cm}
	}
	\end{minipage}%
	\caption{Some results of $\left\|\nabla\cdot \bu_{h}^{k}\right\|_{L^2}$, 
		$\tilde{\mathcal{E}}_{h}^k$, 
		$\left|\dfrac{\tilde{\mathcal{E}}^k_h-\tilde{\mathcal{E}}^{k-1}_h}{\Delta t} - \mathcal{D}^{k+\frac{1}{2}}_h\right|$ over time 
		for the conservation and dissipation tests at $\mathsf{c} = 1$, $N=2$, $h=\frac{1}{8}$ and $\Delta t=\frac{1}{50}$, where 
		$\mathcal{D}_{h}^{k-\frac{1}{2}} := - \Rn^{-1}\mathcal{S}_{h}^{k-\frac{1}{2}}
		- \mathsf{c}\Rm^{-1}\tilde{\mathcal{J}}_{h}^{k-\frac{1}{2}}
		+ \mathsf{c}\left( \mathcal{A}_{h}^{k-\frac{1}{2}} -\tilde{\mathcal{A}}^{k-\frac{1}{2}}_{h}\right)$, cf. \eqref{Eq: discrete energy dissipation}. \revYi{The discrete energy $\mathcal{E}^k_h + (\mathcal{E}^k_h-\mathcal{E}^0_h)\times10^{10}$  ($\color{blue}---$) for the ideal case of a Crank-Nicolson discretization of the coupled formulation (see \ref{App: coupled discretization}) is also shown in the top part of the bottom-left diagram.}
		}
	\label{fig:Energy tests}
\end{figure}

We test the conservation and dissipation properties of the decoupled formulation in a spatial domain $\Omega=\left(x, y, z\right)\in \left[0, 1\right]^3$ using manufactured solutions. The initial conditions, taken from \cite[Section~4.2]{HU2021110284}, are 
\[
\bu^{0}=\begin{bmatrix}
	-\sin\left(\pi(x-\frac{1}{2})\right)\cos\left(\pi(y-\frac{1}{2})\right)z\left(z-1\right) \\ \cos\left(\pi(x-\frac{1}{2})\right)\sin\left(\pi(y-\frac{1}{2})\right)z\left(z-1\right) \\ 0
\end{bmatrix}
\]
and
\[
\bH^{0}=\begin{bmatrix}
	-\sin\left(\pi x\right)\cos\left(\pi y\right) & \cos\left(\pi x\right)\sin\left(\pi y\right)  & 0
\end{bmatrix}^{\mathsf{T}}
\]
which possess the initial energy of $\mathcal{E}^{0} \approx 0.25833$ if $\mathsf{c} = 1$. The boundary conditions are 
\[
\widehat{\bu}=\widehat{\boldsymbol{E}}=\boldsymbol{0}\quad\text{and}\quad \widehat{P} = 0 \quad \text{on } \partial\Omega\times(0,T];
\]
we use natural boundary conditions on the whole boundary and set them to zero for these tests. And the external body force is zero; $\bF=\boldsymbol{0}$. A mesh of uniform cubic elements and a constant time-step interval is employed. The edge length of the element is denoted by $h$. The decoupled formulation, \eqref{Eq: wf step 1} and \eqref{Eq: wf step 2}, is solved with different $\Rn$ and $\Rm$. The results presented in Fig.~\ref{fig:Energy tests} show that mass conservation is always strongly satisfied. 
As for the energy, a dissipation rate of a machine precision error is obtained.

As can be seen in Figure~\ref{fig:Energy tests}, the total energy in the absence of dissipative terms, viscosity, and resistivity is not constant, although it does not deviate too far from its initial value. This is due to the fact that $\mathcal{A}_{h}^{k-\frac{1}{2}} -\tilde{\mathcal{A}}^{k-\frac{1}{2}}_{h} \neq 0$ as was shown in Section~\ref{Sec: discretization}. However, the variations in total energy remain bounded, 
as is expected from a symplectic integrator.

\subsection{Orszag-Tang vortex} \label{Subsec: OTV}
The Orszag-Tang vortex initially studied in \cite{Orszag_Tang_1979} is a well-known two-dimensional incompressible MHD test case. In a periodic square, the Orszag-Tang vortex gradually develops narrow corridors (that will eventually become singularities) of extreme current density where the magnetic field changes its sign suddenly \cite{kraus2018variational}, which makes this test case a challenging one. 

\begin{figure}[!htb]
	\centering
	\begin{minipage}[c]{0.92\textwidth}
		\centering{
			\subfloat{
				\begin{minipage}[b]{0.31\textwidth}
					\centering
					\includegraphics[width=0.95\linewidth]{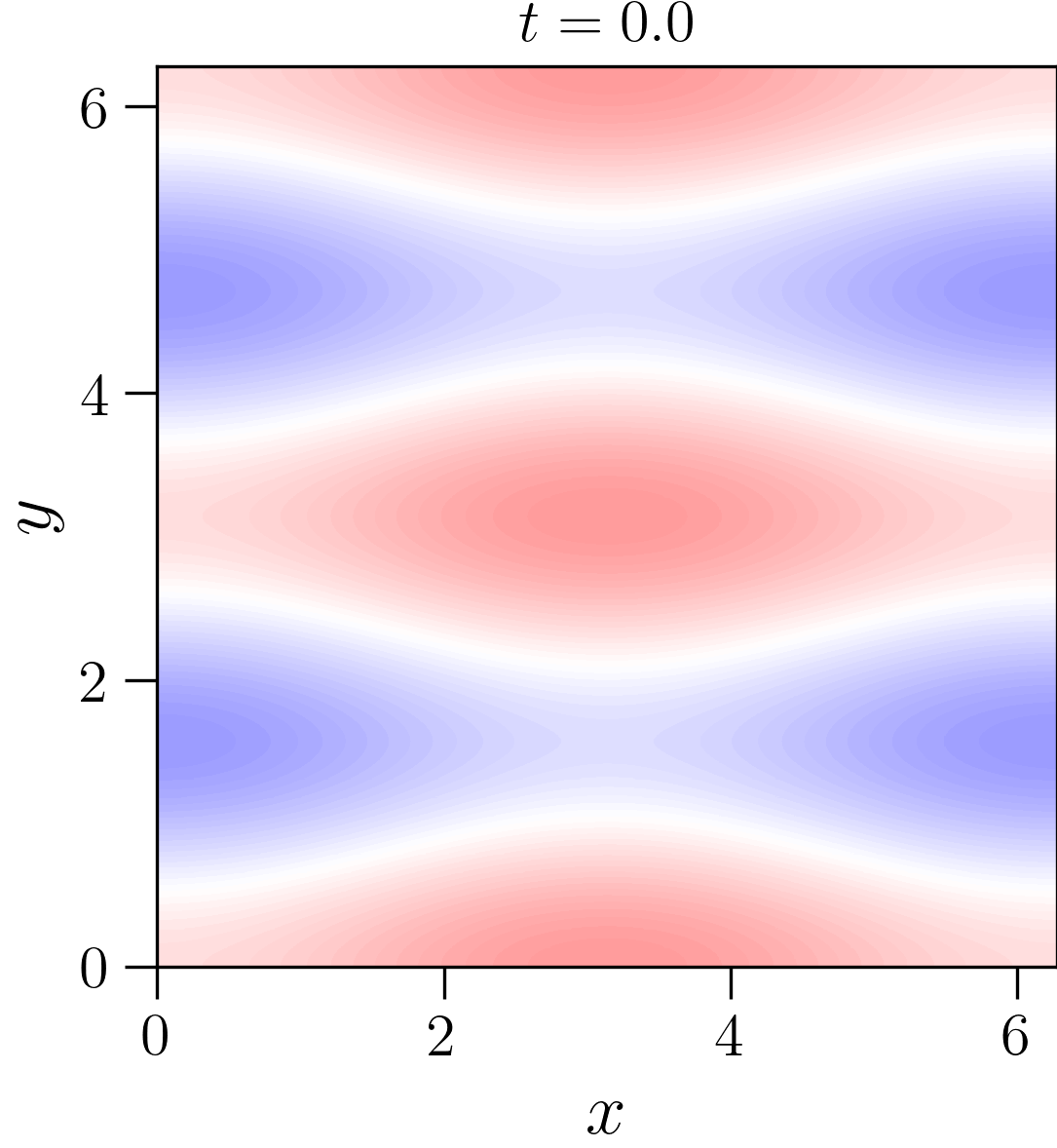}
				\end{minipage}
			}
			\subfloat{
				\begin{minipage}[b]{0.31\textwidth}
					\centering
					\includegraphics[width=0.95\linewidth]{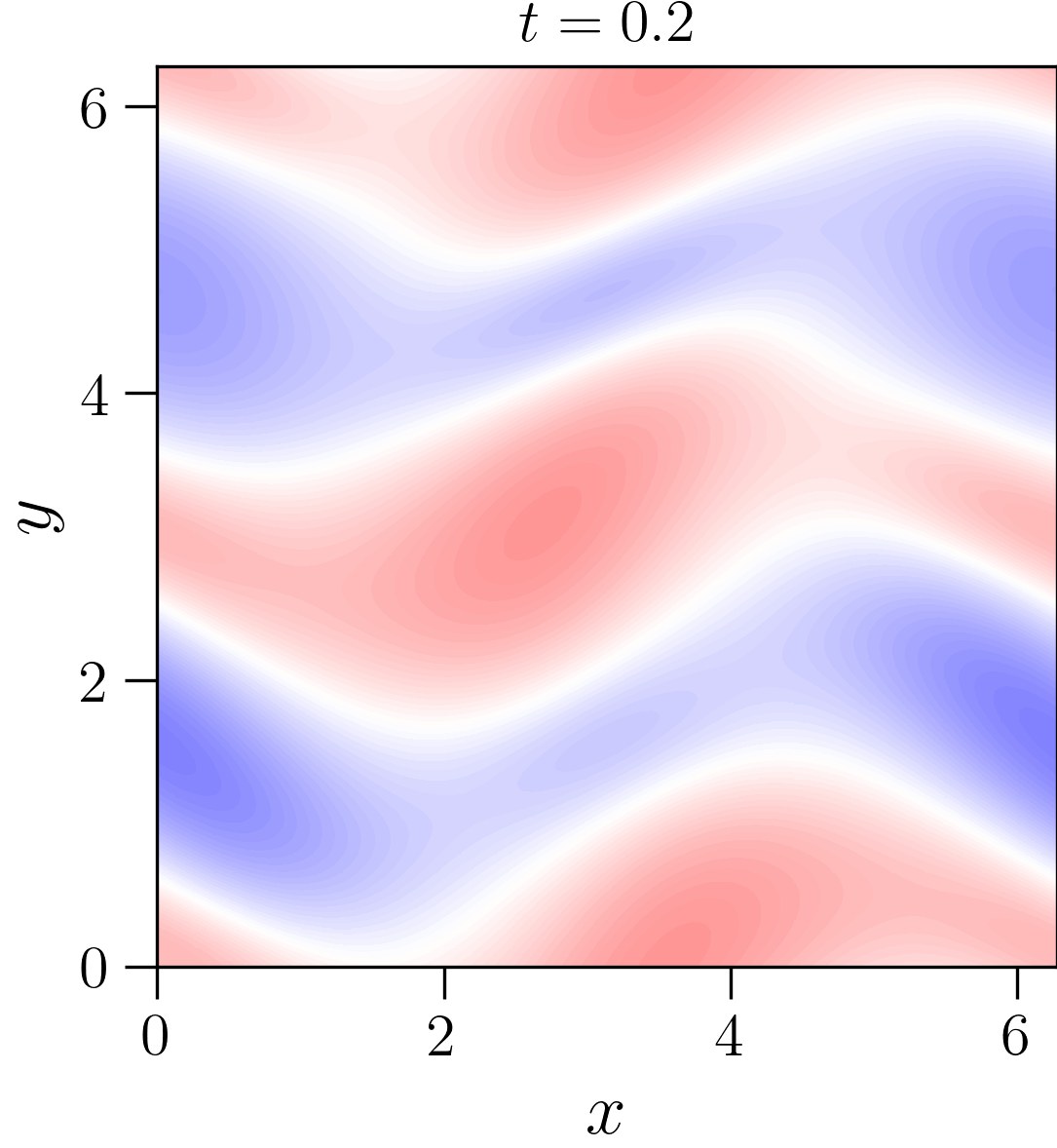}
				\end{minipage}
			}
			\subfloat{
				\begin{minipage}[b]{0.31\textwidth}
					\centering
					\includegraphics[width=0.95\linewidth]{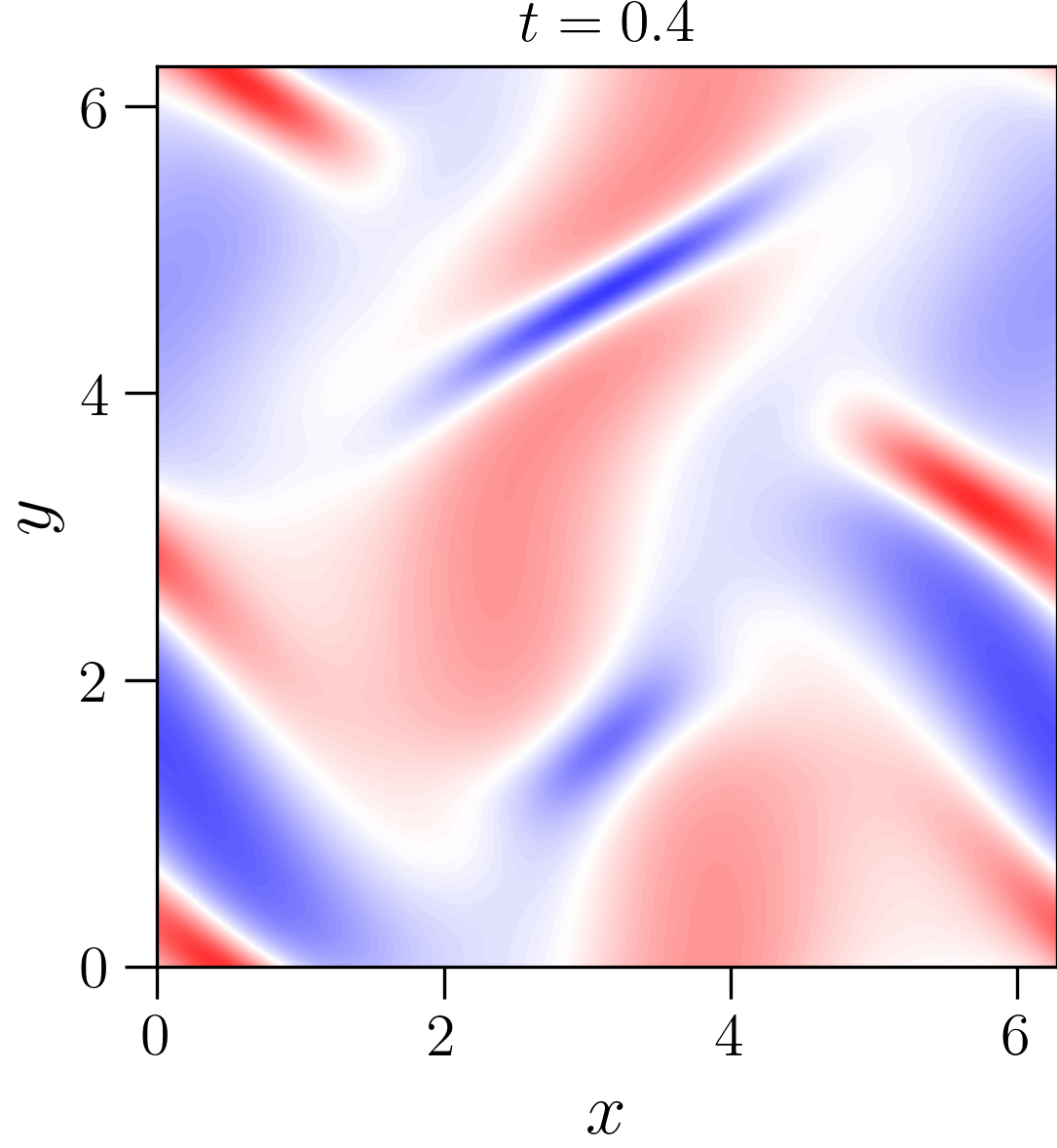}
				\end{minipage}
			}\\
			\subfloat{
				\begin{minipage}[b]{0.31\textwidth}
					\centering
					\includegraphics[width=0.95\linewidth]{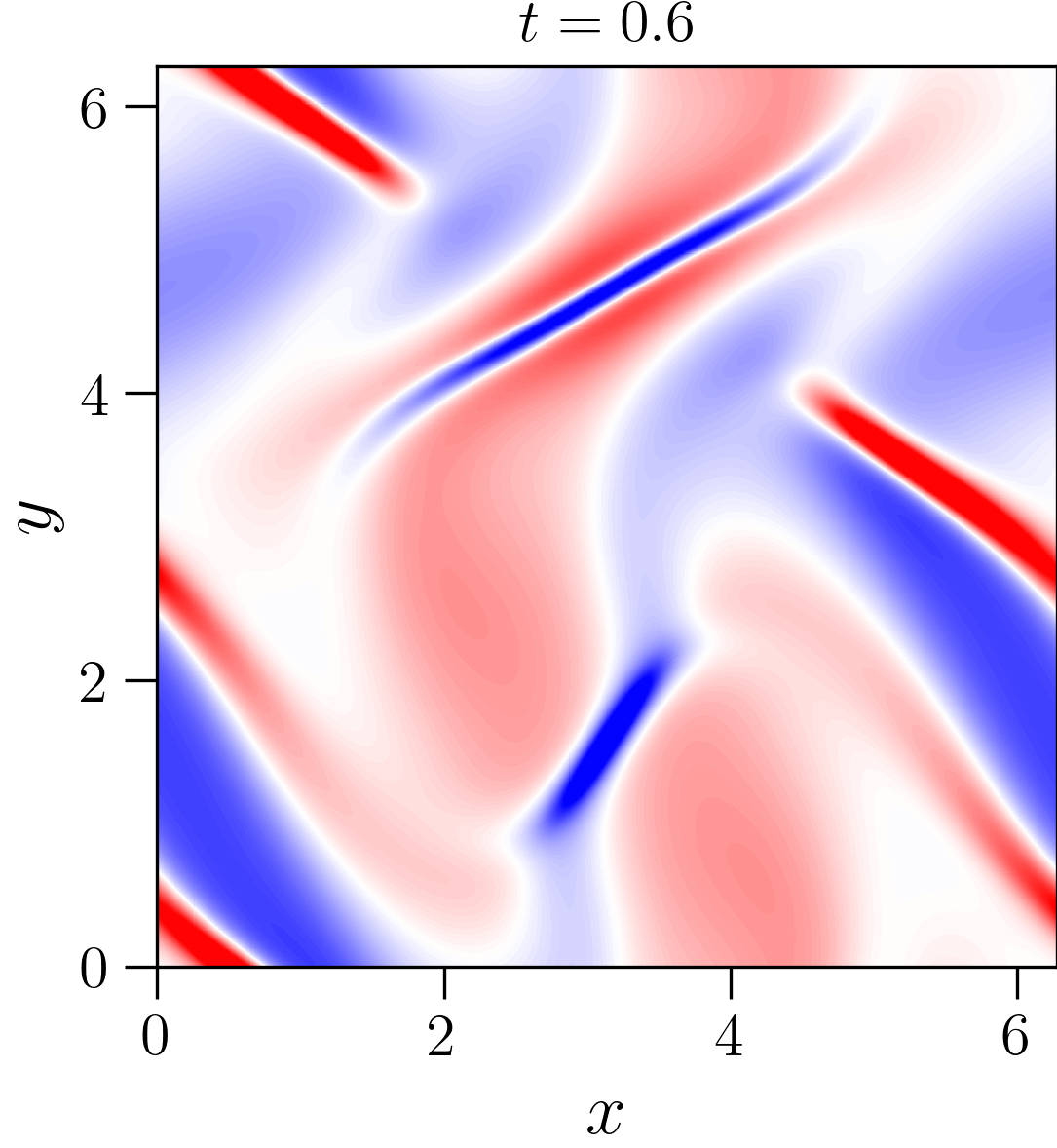}
				\end{minipage}
			}
			\subfloat{
				\begin{minipage}[b]{0.31\textwidth}
					\centering
					\includegraphics[width=0.95\linewidth]{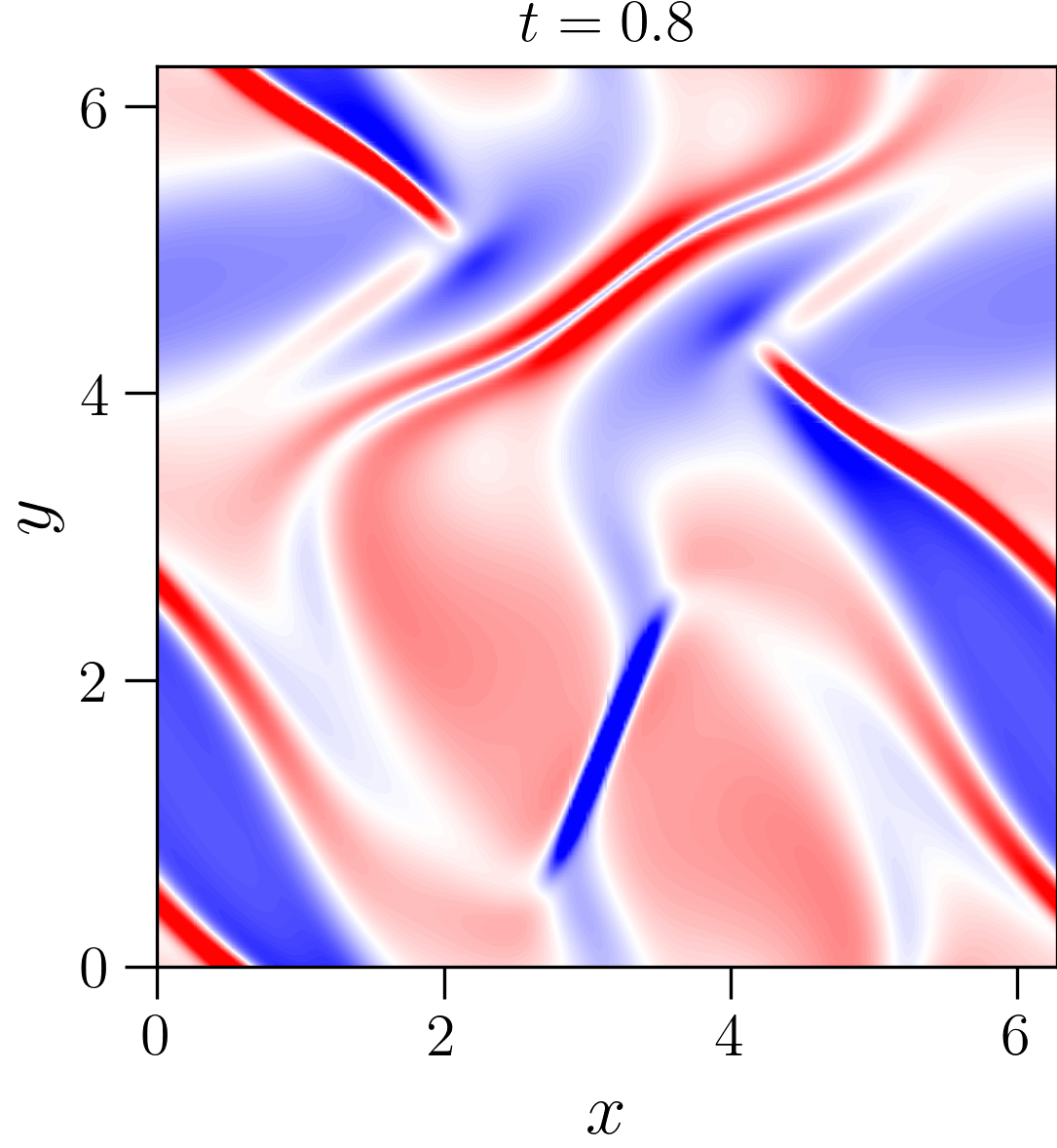}
				\end{minipage}
			}
			\subfloat{
				\begin{minipage}[b]{0.31\textwidth}
					\centering
					\includegraphics[width=0.95\linewidth]{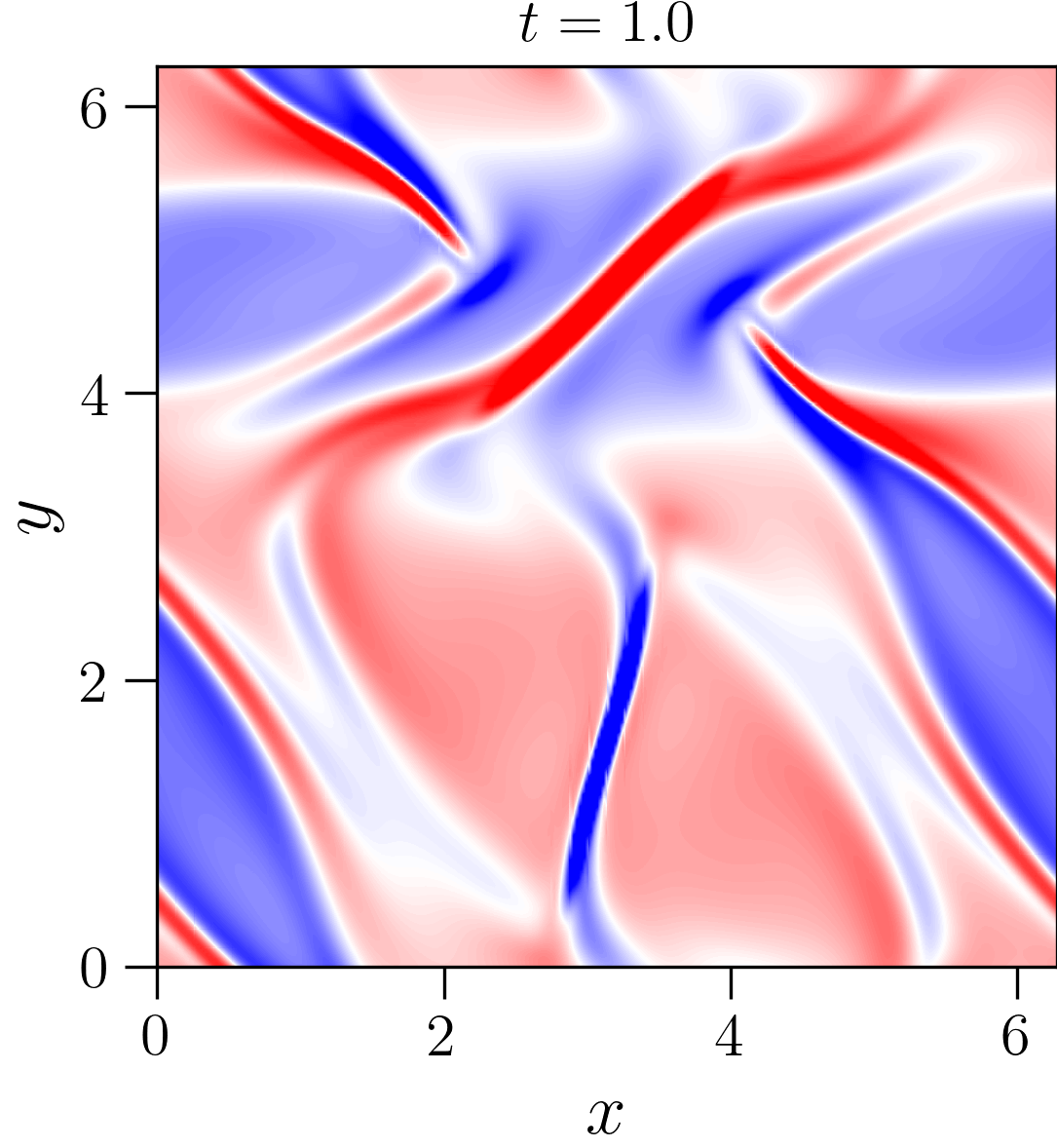}
				\end{minipage}
			}
		}
	\end{minipage}
	\begin{minipage}[c]{0.06\textwidth}%
		\centering
		\includegraphics[width=0.7\linewidth,]{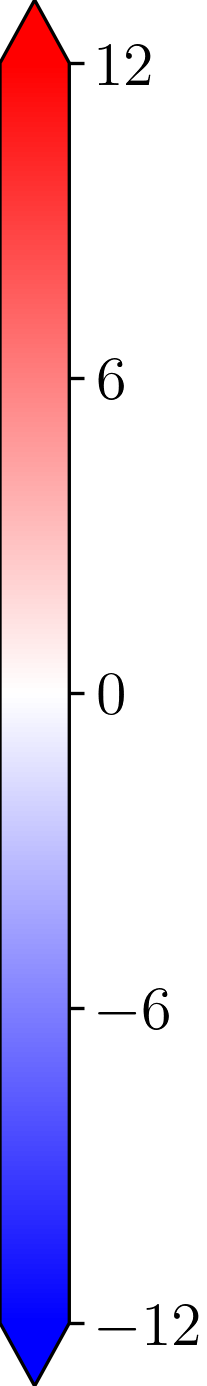}
	\end{minipage}
\caption{$\boldsymbol{j}_{h}^{k+\frac{1}{2}}$ at $t^{k}=t\in \left\lbrace0, 0.2, 0.4, 0.6, 0.8, 1\right\rbrace$ for the Orszag-Tang vortex test using $N=4, K=48, \Delta t=\frac{1}{200}$.}
\label{fig:OTV tests}
\end{figure}

\begin{figure}[!htb]
	\centering
	\begin{minipage}[c]{0.92\textwidth}
		\centering{
			\subfloat{
				\begin{minipage}[b]{0.31\textwidth}
					\centering
					\includegraphics[width=0.95\linewidth]{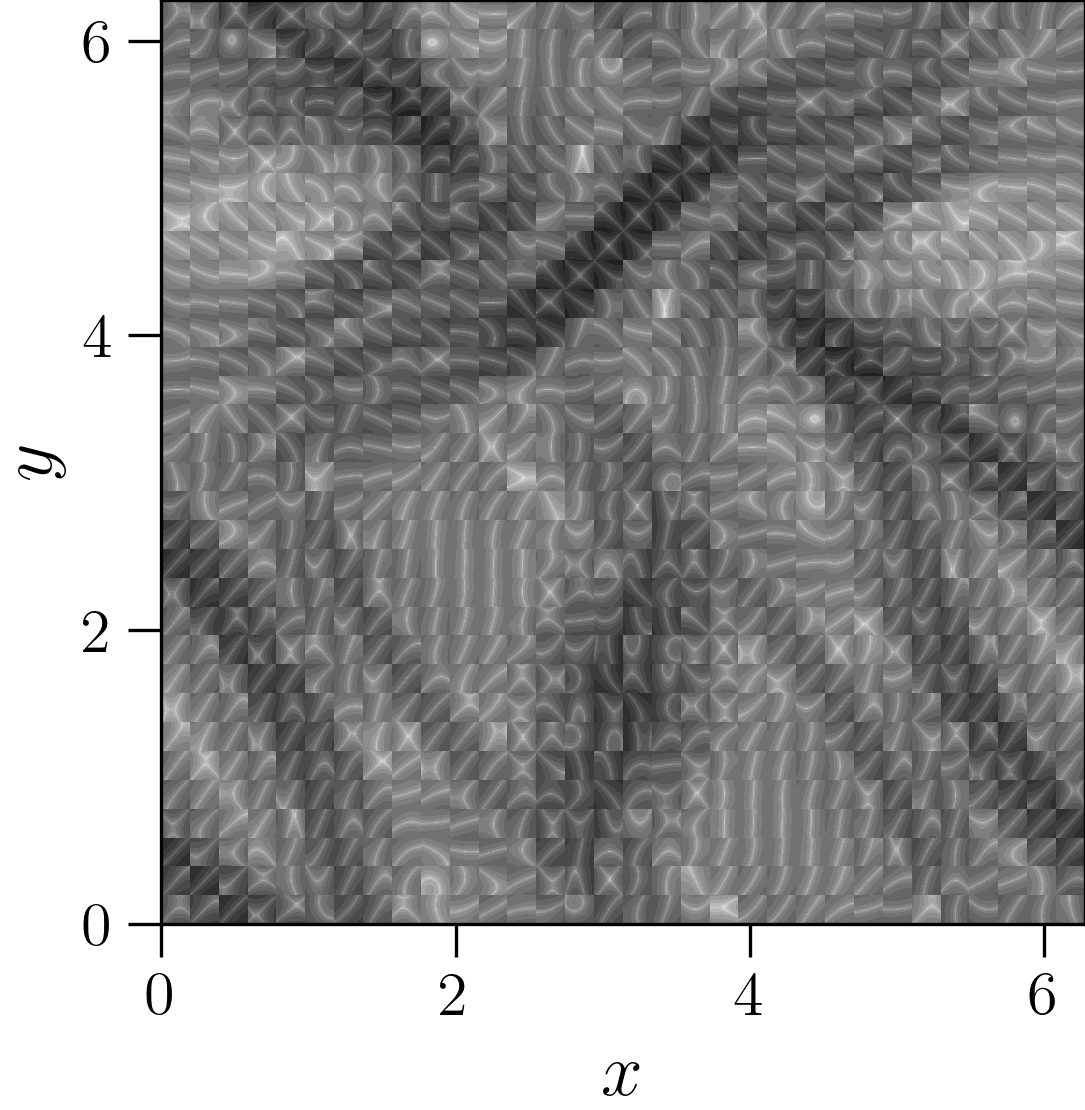}
				\end{minipage}
			}
			\subfloat{
				\begin{minipage}[b]{0.31\textwidth}
					\centering
					\includegraphics[width=0.95\linewidth]{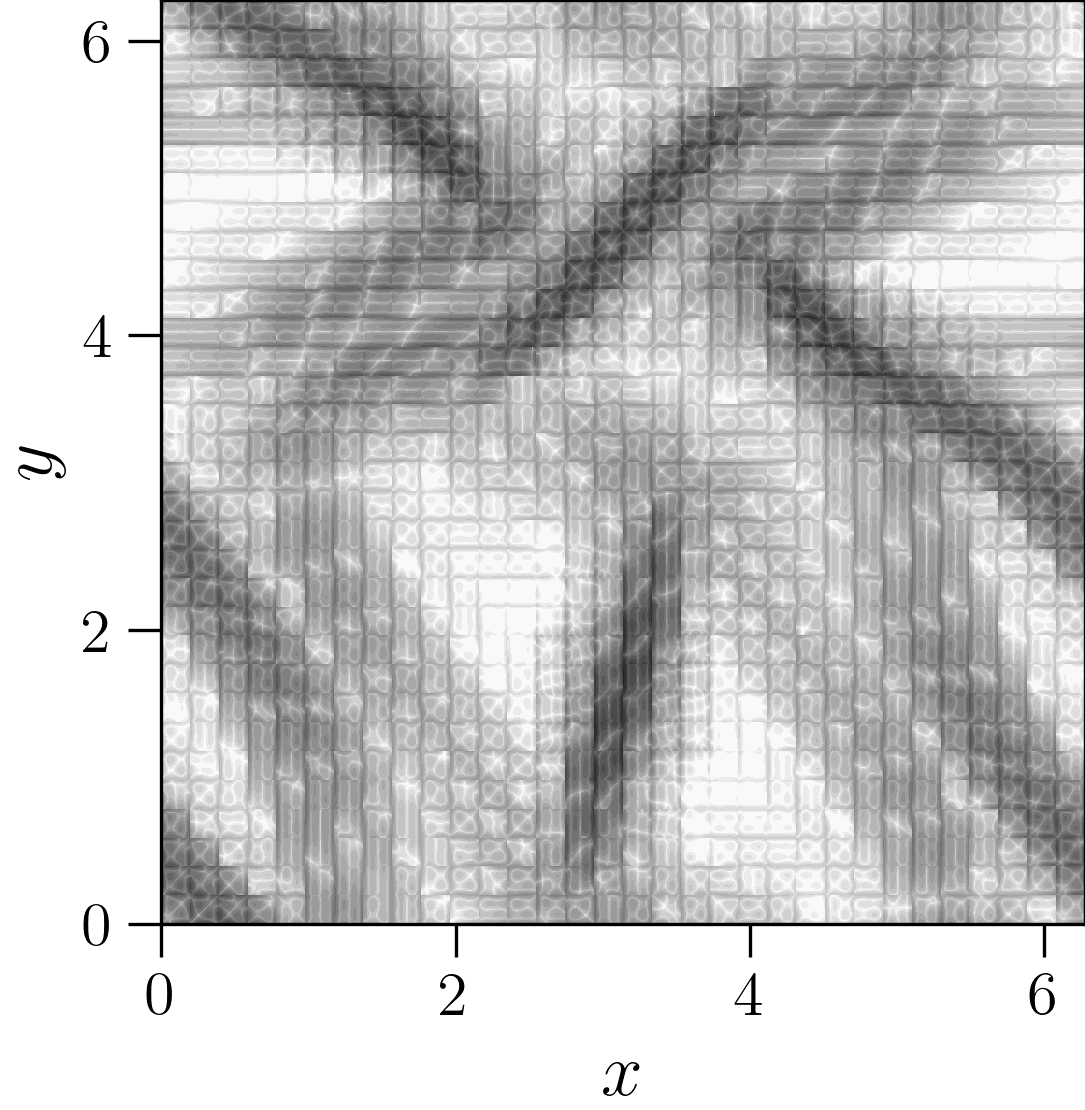}
				\end{minipage}
			}
			\subfloat{
				\begin{minipage}[b]{0.31\textwidth}
					\centering
					\includegraphics[width=0.95\linewidth]{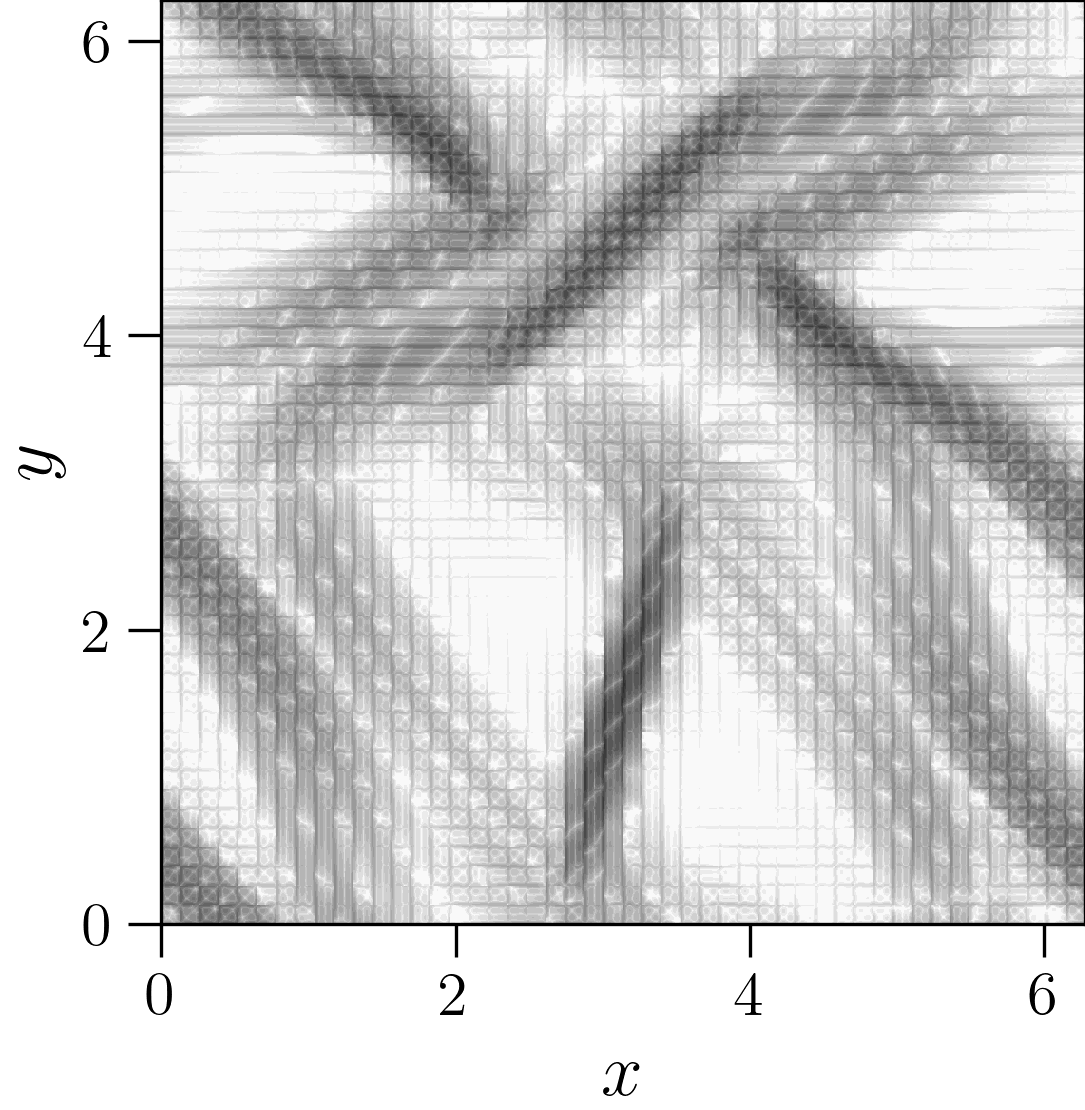}
				\end{minipage}
			}
		}
	\end{minipage}
	\begin{minipage}[c]{0.05\textwidth}%
		\centering
		\includegraphics[width=0.7\linewidth]{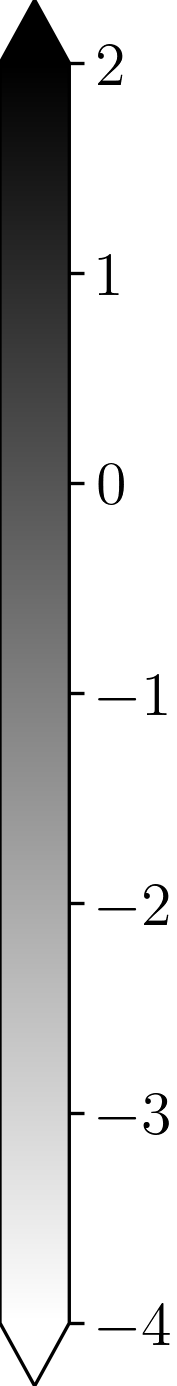}
	\end{minipage}
	\caption{Element-wise $\log_{10}\left(\left|\nabla\cdot\boldsymbol{H}_{h}^{k+\frac{1}{2}}\right|\right)$ at $t^{k}=t=1$ for the Orszag-Tang vortex test using $\Delta t = \frac{1}{200}$. Left: $N=2, K=32$. Middle: $N=4, K=32$. Right: $N=4, K=48$.}
	\label{fig:OTV divH tests}
\end{figure}

Giving a stream function $\psi$ and a magnetic potential $A$, also see \cite[Section 5.3]{kraus2018variational},
\[
\psi = 2\sin\left(y\right) - 2 \cos\left(x\right),\quad A = \cos\left(2y\right) - 2\cos\left(x\right),
\]
the initial conditions are $\bu^{0}=\nabla\times\psi $, 
and $\bH^{0}=\nabla\times A$.
The fully orthogonally periodic square is $\Omega = \left(x, y\right)\in \left[0, 2\pi\right]^2$ and $\bF=\boldsymbol{0}$. On a mesh of $K^2$ uniform cubic elements, we solve the two-dimensional version of the decoupled formulation with a time-step interval $\Delta t=\frac{1}{200}$ and the parameters are $\Rn = 100$, $\Rm=100$, and $\mathsf{c}=1$. 
Some snapshots of $\boldsymbol{j}_{h}^{k+\frac{1}{2}} = \nabla\times \bH _{h}^{k+\frac{1}{2}}$ are shown in Fig.~\ref{fig:OTV tests}, and a good match to reference results in \cite{kraus2018variational} is observed. Elementwise $\log_{10}\left(\left|\nabla\cdot\boldsymbol{H}_{h}^{k+\frac{1}{2}}\right|\right)$ at $t^k=t=1$ for different combinations of $N$ and $K$ are shown in Fig.~\ref{fig:OTV divH tests} where we can see the convergence of weak conservation of Gauss's law for magnetism under $ph$-refinement.

\subsection{Lid-driven cavity} \label{Subsec: LDC}
The two-dimensional lid-driven cavity is a domain $ \Omega = \left(x, y\right)\in [0,1]^2 $ with four infinite no-slip solid walls, and there is no external body force; $\bF=\boldsymbol{0}$. Initially, the fluid is at rest, i.e., $\bu^{0}=\boldsymbol{0}$, and is immersed in a magnetic field $\bH^{0} = \begin{bmatrix}
0 & 1
\end{bmatrix}^{\mathsf{T}}$. At $t>0$, the top wall, namely the lid, moves right with a constant speed $1$, i.e., $\bu_{\mathrm{lid}} = \begin{bmatrix}
	1 & 0
\end{bmatrix}^{\mathsf{T}}$, and drives the fluid due to viscous entrainment. Over the whole boundary, the tangential electric field is zero all the time. In summary, we have boundary conditions
\[
	\bu\cdot\bn = \hat{u} = 0,\quad 
	\bE\times \bn = \widehat{E} = 0
\qquad \text{on } \partial\Omega\times\left(0,T\right],
\]
and 
\[
\bu\times \bn = \hat{u}_{\parallel} = \left\lbrace\begin{aligned}
	1 \quad &\text{on}\ \Gamma_{y}^{+}\times\left(0,T\right] \\
	0 \quad&\text{else}
\end{aligned}\right.,
\]
where $\Gamma_{y}^{+}$ is the face $\left(x, y\right)\in \left(0,1\right)\times1$. 

\begin{figure}[h!]
	\centering
	\begin{minipage}[c]{1\textwidth}
		\centering{
			\subfloat[mesh]{ \label{LDC a}
				\begin{minipage}[b]{0.5\textwidth}
					\centering
					\includegraphics[width=0.68\linewidth]{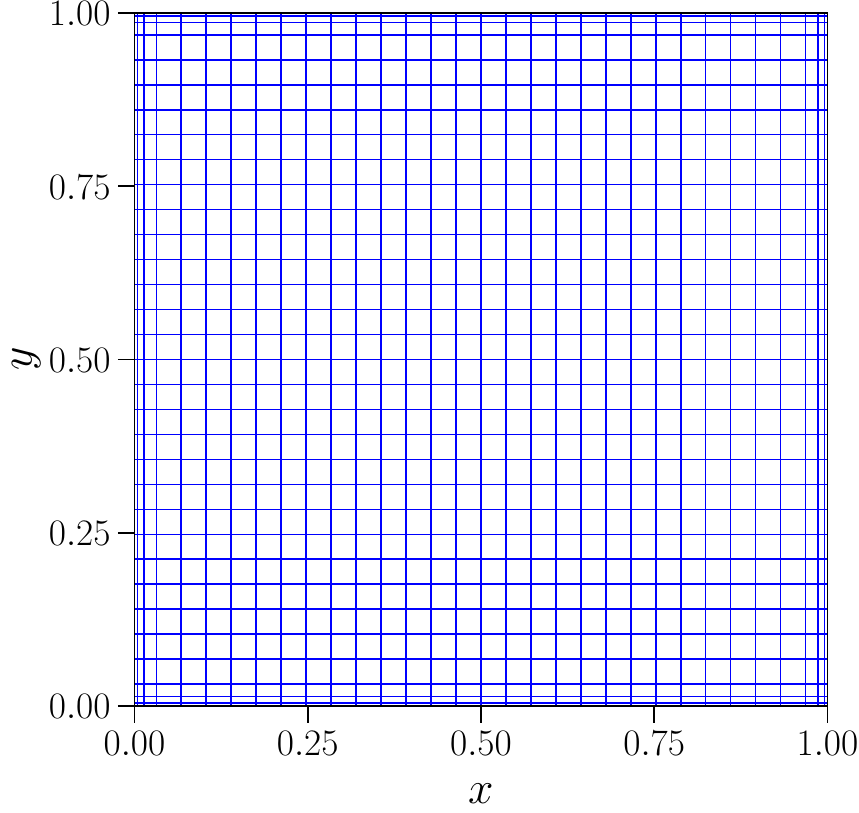}
				\end{minipage}
			}
			\subfloat[$\omega_{h}^k$]{ \label{LDC b}
				\begin{minipage}[b]{0.5\textwidth}
					\centering
					\includegraphics[width=0.68\linewidth]{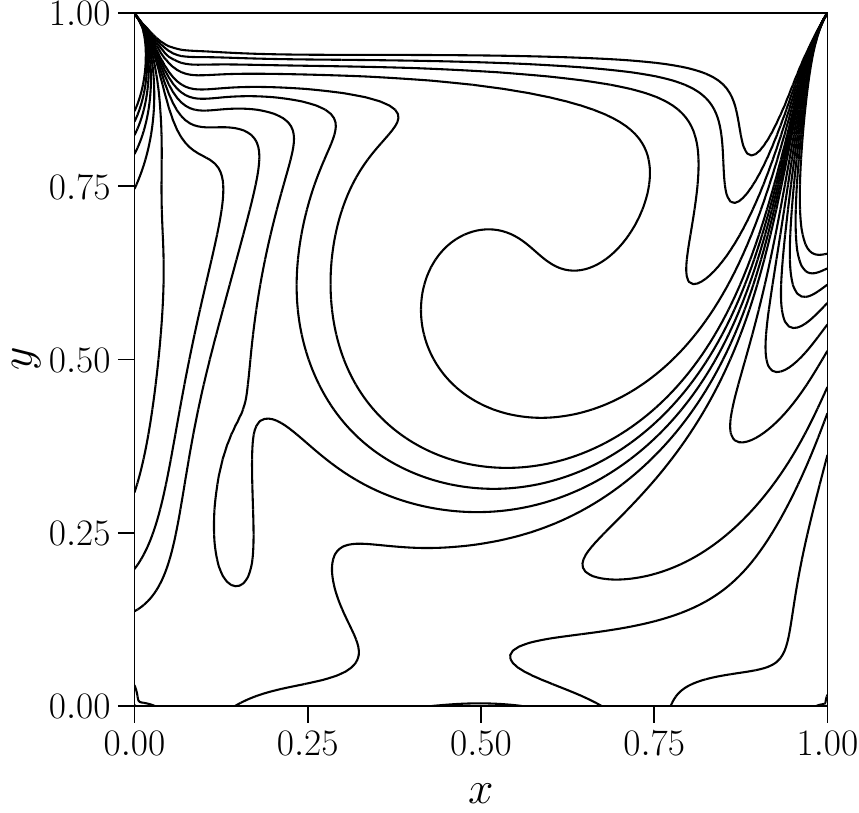}
				\end{minipage}
			}\\
			\subfloat[streamlines of $\bu_{h}^{k}$]{ \label{LDC c}
				\begin{minipage}[b]{0.5\textwidth}
					\centering
					\includegraphics[width=0.68\linewidth]{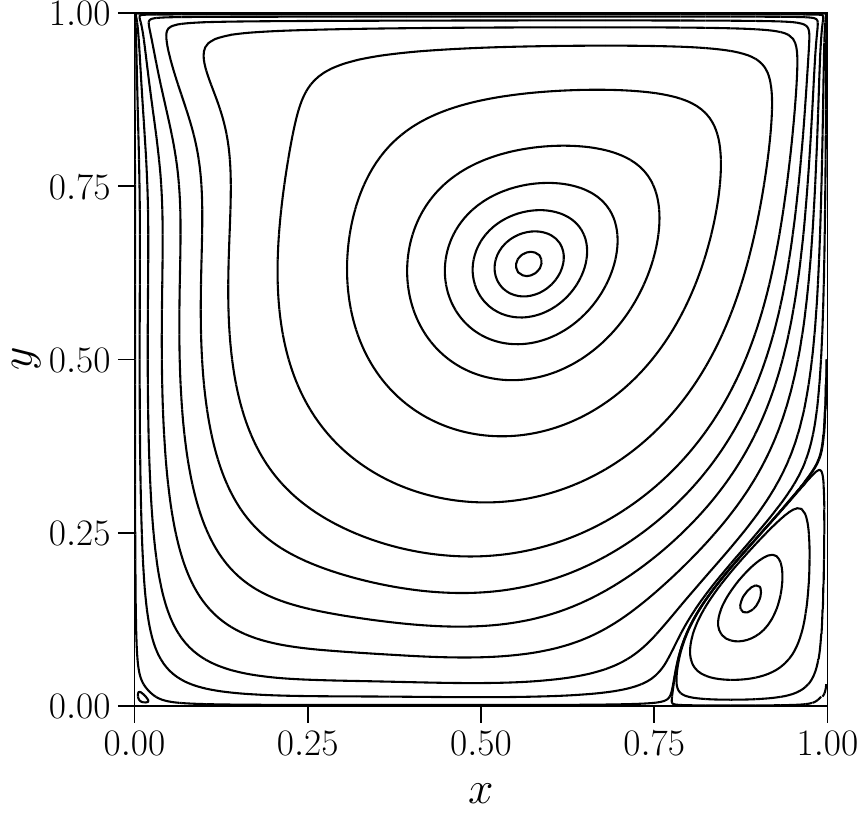}
				\end{minipage}
			}
			\subfloat[streamlines of $\bH_{h}^{k+\frac{1}{2}}$]{ \label{LDC d}
				\begin{minipage}[b]{0.5\textwidth}
					\centering
					\includegraphics[width=0.68\linewidth]{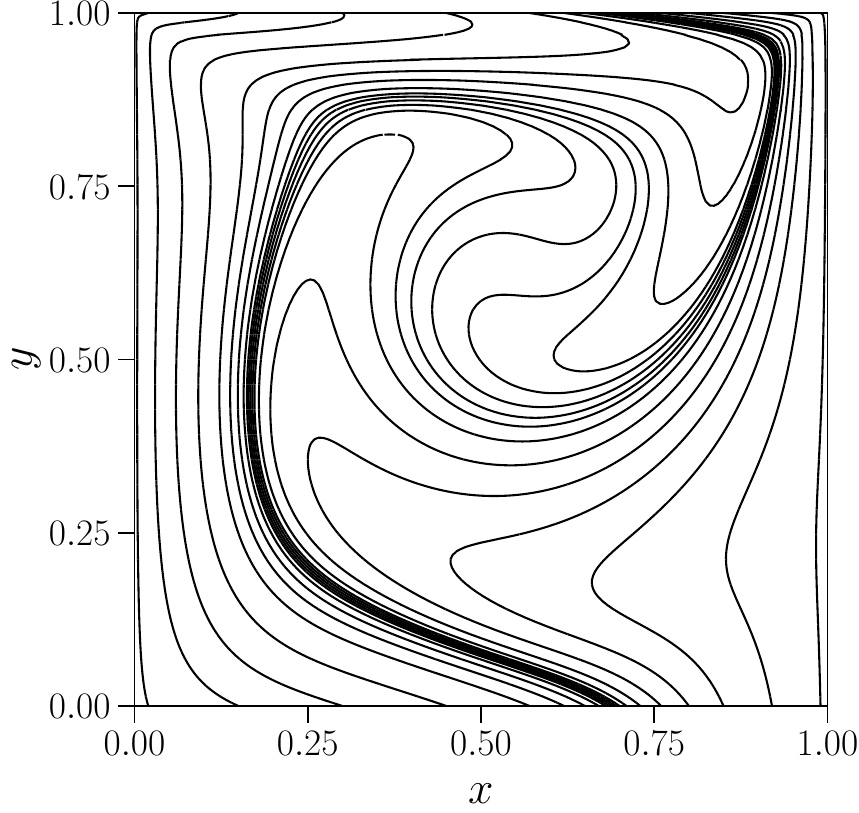}
				\end{minipage}
			}
		}
	\end{minipage}%
	\caption{The mesh and some results at the steady state for the lid-driven cavity test. (a): The mesh. (b): Vorticity $\omega_{h}^{k}$ for contour lines in $\left\lbrace -5, -4, -3, -2, -1, -0.5, 0, 0.5, 1, 2, 3, 4, 5, 6, 7\right\rbrace$. (c): Streamlines of $\bu_{h}^k$. (d): Streamlines of $\bH_{h}^{k+\frac{1}{2}}$. Note that, for diagram (d), a direct numerical integral over a grid of discrete values of $\bH_{h}^{k+\frac{1}{2}}$ is applied. Since $\bH_{h}^{k+\frac{1}{2}}$ is not strongly divergence free, there does not exist a stream function $\psi$ such that $\boldsymbol{H}_h^{k+\frac{1}{2}}=\nabla\times\psi$. However, as we have shown that $\bH_{h}^{k+\frac{1}{2}}$ is weakly divergence free for the proposed method, see Fig.~\ref{fig:divH convergence rate} and Fig.~\ref{fig:OTV divH tests}, the error in the stream function by this approach is limited. Stream functions of contour lines in (c) and (d) are given in See \ref{App: streamfunctions}. A good match to reference results in \cite[Fig.~6, Fig.~7]{ZHANG202345} is obtained.
	} 
	\label{fig:LDC}
\end{figure}

To simulate the lid-driven cavity, an orthogonal spatial mesh of $32^2$ elements that are locally refined near the boundary, see Fig.~\ref{LDC a}, is used. The time-step interval is $\Delta t = \frac{1}{1000}$, Reynolds numbers are $\Rn=\Rm=400$, the coupling number is $\mathsf{c}= \Rm^{-1}$ (such that $\mathsf{s} = \mathsf{c}\Rm = 1$)
and the polynomial degree is $N=3$. The two-dimensional version of the decoupled formulation is solved until a steady state is reached. The criterion of the steady state is 
\[
\frac{1}{\Delta t}\max\left(\left\|\bu_{h}^{k}-\bu_{h}^{k-1}\right\|_{L^2}, \left\|\bH_{h}^{k+\frac{1}{2}}-\bH_{h}^{k-\frac{1}{2}}\right\|_{L^2}\right) < 10^{-5}.
\]
The flow reaches its steady state at $t=82.049$. Some results are presented in Fig.~\ref{fig:LDC}. They have a good agreement to the reference results taken from \cite{ZHANG202345}. In Table~\ref{Tab: LDC centerline x} and Table~\ref{Tab: LDC centerline y}, solutions along two centerlines are provided as quantitative benchmarks.

\begin{table}[h!]
	\caption{Solutions along the $x$-direction centerline, i.e., $\left(x, y\right) = [0, 1]\times 0.5$, at the steady state for the lid-driven cavity test. $\boldsymbol{u}_{h}^{k} = \begin{bmatrix}
		u_{h}^{k} & v_{h}^{k}
	\end{bmatrix}^\mathsf{T}$ and $\bH _{h}^{k+\frac{1}{2}} = \begin{bmatrix}
	H_{x; h}^{k+\frac{1}{2}} & H_{y;h}^{k+\frac{1}{2}}
	\end{bmatrix}^\mathsf{T}$.}
	\centering
	\begin{tabular}{lrrrrrrr}
		\hline\bigstrut
		$x$  & &\multicolumn{1}{r}{$u^k_h$}&\multicolumn{1}{r}{$v^k_h$} & \multicolumn{1}{r}{$\omega^k_h$} & \multicolumn{1}{r}{$H^{k+\frac{1}{2}}_{x;h}$}& \multicolumn{1}{r}{$H^{k+\frac{1}{2}}_{y;h}$} \\ \hline\bigstrut
		0    &&   $ 0 $     &    $ 0 $ &   $ 2.88484$  & $ 0 $     &$ 4.98952$ \\
		0.05 &&  $-0.00152$ & $0.11308$&   $ 1.63098$  & $0.09510$ &$ 4.90996$ \\
		0.1  &&  $-0.01043$ & $0.17267$&   $ 0.63099$  & $0.15636$ &$ 4.15657$ \\
		0.15 &&  $-0.02841$ & $0.19891$&    $0.09245$  & $0.14927$ &$ 2.49169$ \\
		0.25 &&  $-0.07561$ & $0.22147$&   $-0.48799$  & $0.13626$ &$-0.07002$ \\
		0.5  &&  $-0.14277$ & $0.05402$&   $-2.19175$  & $0.05387$ &$-0.07419$ \\
		0.75 &&  $-0.20138$ &$-0.25127$&   $-2.33726$  & $0.41737$ &$ 0.35175$ \\
		0.85 &&  $-0.16772$ &$-0.33732$&   $ 0.53077$  & $0.73888$ &$ 1.37109$ \\
		0.9  &&  $-0.09707$ &$-0.26638$&   $ 2.68539$  & $0.56165$ &$ 1.55770$ \\
		0.95 &&  $-0.02697$ &$-0.12420$&   $ 3.14441$  & $0.27997$ &$ 1.49902$ \\
		1    && 	$ 0 $   &  $ 0 $   &   $ 1.72384$  & $ 0 $     &$ 1.46331$ \\ \hline
	\end{tabular} 
	\label{Tab: LDC centerline x}
\end{table}

\begin{table}[h!]
	\caption{Solutions along the $y$-direction centerline, i.e., $\left(x, y\right) = 0.5\times [0, 1]$, at the steady state for the lid-driven cavity test. $\boldsymbol{u}_{h}^{k} = \begin{bmatrix}
			u_{h}^{k} & v_{h}^{k}
		\end{bmatrix}^\mathsf{T}$ and $\bH _{h}^{k+\frac{1}{2}} = \begin{bmatrix}
			H_{x; h}^{k+\frac{1}{2}} & H_{y;h}^{k+\frac{1}{2}}
		\end{bmatrix}^\mathsf{T}$.}
	\centering
	\begin{tabular}{lrrrrrr}
		\hline\bigstrut
		$y$  & &\multicolumn{1}{r}{$u^k_h$}&\multicolumn{1}{r}{$v^k_h$} & \multicolumn{1}{r}{$\omega^k_h$}  & \multicolumn{1}{r}{$H^{k+\frac{1}{2}}_{x;h}$}& \multicolumn{1}{r}{$H^{k+\frac{1}{2}}_{y;h}$} \\ \hline\bigstrut
		0    &&       $ 0 $ &    $ 0 $  &$1.06180$  & $-2.74879$ &    $ 1 $  \\
		0.05 &&   $-0.03939$& $-0.00056$&$0.60008$  & $-2.40069$ &$ 0.80054$ \\
		0.1  &&   $-0.07135$& $-0.00329$&$0.58025$  & $-1.65057$ &$ 0.60078$ \\
		0.15 &&   $-0.11225$& $-0.00676$&$0.70413$  & $-0.73480$ &$ 0.39444$ \\
		0.25 &&   $-0.21581$& $-0.00329$&$ 0.35426$  & $ 0.52230$ &$ 0.06335$ \\
		0.5  &&   $-0.14277$& $ 0.05402$&$-2.19175$  & $ 0.05387$ &$-0.07419$ \\
		0.75 &&   $ 0.12722$& $ 0.07170$&$-1.81455$  & $-0.11813$ &$-0.04774$ \\
		0.85 &&   $ 0.22351$& $ 0.05999$&$-1.41003$  & $ 0.33644$ &$-0.08662$ \\
		0.9  &&   $ 0.28772$& $ 0.04373$&$-2.18569$  & $ 1.71570$ &$-0.08702$ \\
		0.95 &&   $ 0.47275$& $ 0.01801$&$-6.41153$  & $ 3.05614$ &$ 0.20349$ \\
		1    &&      $ 1 $  &    $ 0 $  &$-15.51828$ & $-5.29567$ &   $ 1 $   \\ \hline
	\end{tabular}
	\label{Tab: LDC centerline y}
\end{table}

\section{Conclusion} \label{Sec: conclusion}
In this work, we present a decoupled structure-preserving discretization for the incompressible MHD equations with general boundary conditions. A spatially discrete formulation in a mixed finite element setting is first proposed. It preserves conservation of mass and conservation of charge strongly, preserves Gauss's law for magnetism weakly, and preserves the correct energy dissipation rate. A leapfrog-type temporal integrator is then applied to the fluid part and the Maxwell part at two staggered time sequences, respectively, such that they are decoupled at the fully discrete level. This two-step decoupled method has an optimal spatial accuracy and a second-order temporal accuracy. 
However, the decoupled method does not preserve energy conservation. Other invariants, for example, cross- and magnetic-helicity, are preserved by the incompressible MHD equations. Possible extensions of this work include the preservation of these quantities.

\section*{Acknowledgments}
The research of Deepesh Toshniwal is
supported by project number 212.150 awarded through the Veni research programme by the Dutch Research Council (NWO). Andrea Brugnoli acknowledges the financial support of the German Research Foundation (DFG) and the Berlin Mathematics Research Center (MATH+ project AA4-12). The research of Yi Zhang is supported by the Natural Science Foundation of Guangxi under grant number 2024JJB110005.

\appendix
\section{Stream functions} \label{App: streamfunctions}
In Fig.~\ref{LDC c}, the streamlines are for stream functions in
\[
\begin{aligned}
	\left\lbrace-0.093, -0.092, -0.09, -0.086, -0.078, -0.06, -0.035, -0.018, -0.01, -5\text{E}-3, \right.\qquad \\
	\left.-2\text{E}-3, -5\text{E}-4, -8\text{E}-5, -1\text{E}-6, 5\text{E}-8, 6\text{E}-6, 8\text{E}-5, 3\text{E}-4, 4.1\text{E}-4 \right\rbrace.
\end{aligned}
\]
In Fig.~\ref{LDC d}, the streamlines are for stream functions in
\[
\begin{aligned}
	\left\lbrace-0.99, -0.92, -0.85, -0.80, -0.76, -0.73,-0.71, -0.7, -0.695, -0.69, \right.\qquad \\
	\left.-0.685, -0.68, -0.67, -0.65, -0.62, -0.57, -0.45, -0.3, -0.15, -0.02 \right\rbrace.
\end{aligned}
\]
For both cases, the reference stream function is set to $0$ at the bottom-left corner.

\revYi{
\section{A Crank-Nicolson discretization of the coupled formulation} \label{App: coupled discretization}
On the time sequence \eqref{eq: time sequence}, if we apply the Crank-Nicolson scheme directly to the coupled formulation \eqref{Eq: wf}, we will obtain a fully discrete formulation written as:  Given $\bF\in\left[L^2(\Omega)\right]^3$, natural boundary conditions $\widehat{P}\in H^{1/2}\left(\Omega;\Gamma_{\widehat{P}}\right),\  \widehat{\bu}\in\mathcal{T}H\left(\mathrm{curl};\Omega,\Gamma_{\widehat{\bu}}\right),\  \widehat{\bE}\in \mathcal{T}H\left(\mathrm{curl};\Omega,\Gamma_{\widehat{\boldsymbol{E}}}\right)$,
and initial conditions 
$\left(\bu_{h}^{0},\bH^{0}_{h}\right)\in D\left(\Omega\right)\times C\left(\Omega\right)$, for $k=1,2,3,\cdots$ successively, 
seek $\left(\boldsymbol{u}^k_{h}, \bw^k_{h}, P^{k-\frac{1}{2}}_{h}, \bH^{k}_{h}\right)\in D_{\hat{u}}(\Omega,\Gamma_{\hat{u}})\times C^{\parallel}_{\widehat{\bw}}(\Omega,\Gamma_{\widehat{\bw}})\times S(\Omega)\times C^{\parallel}_{\widehat{\boldsymbol{H}}}(\Omega,\Gamma_{\widehat{\boldsymbol{H}}})$, such that $\forall \left(\bv_{h},\boldsymbol{w}_{h}, q_{h},\boldsymbol{b}_{h}\right)\in D_{0}\left(\Omega,\Gamma_{\hat{u}}\right) \times C^{\parallel}_{\boldsymbol{0}}\left(\Omega,\Gamma_{\widehat{\bw}}\right)\times S\left(\Omega\right)\times C_{\boldsymbol{0}}^{\parallel}\left(\Omega,\Gamma_{\widehat{\bH}}\right) $,
\begin{subequations}
	\begin{align*}
		\left\langle\dfrac{\bu_h^k-\bu^{k-1}_h}{\Delta t}, \bv_{h}\right\rangle_{\Omega}
		+ \mathcal{F}\left(\dfrac{\bw_h^k+\bw_h^{k-1}}{2},\dfrac{\bu_h^k+\bu_h^{k-1}}{2},\dfrac{\bH_h^k+\bH_h^{k-1}}{2},\bv_h\right)\hspace{-1cm}&\\
		- \ip{P^{k-\frac{1}{2}}_{h}}{\nabla\cdot\bv_{h}}{\Omega}  &= \ip{\bF^{k-\frac{1}{2}}}{\bv_{h}}{\Omega} - \ip{\widehat{P}^{k-\frac{1}{2}}}{\mathcal{T}\bv_{h}}{\Gamma_{\widehat{P}}},\\
		-\ip{\bu^k_{h}}{\nabla\times\boldsymbol{w}_{h}}{\Omega}+\ip{\bw^k_{h}}{\boldsymbol{w}_{h}}{\Omega} & = - \ip{\widehat{\bu}^k}{\mathcal{T}_{\parallel}\boldsymbol{w}_{h}}{\Gamma_{\widehat{\bu}}} ,\\
		\ip{\nabla\cdot\bu^k_{h}}{q_{h}}{\Omega} & = 0,\\
		\ip{\dfrac{\bH_h^k-\bH^{k-1}_h}{\Delta t}}{\boldsymbol{b}_{h}}{\Omega} + \mathcal{M}\left(\dfrac{\bw_h^k+\bw_h^{k-1}}{2}, \dfrac{\bH_h^k+\bH_h^{k-1}}{2},\boldsymbol{b}_{h}\right)
		&= \ip{\widehat{\bE}^{k-\frac{1}{2}}}{\mathcal{T}_{\parallel}\boldsymbol{b}_{h}}{\Gamma_{\widehat{\bE}}}.
	\end{align*}
\end{subequations}
where $\mathcal{F}$ and $\mathcal{M}$ are functionals
\[\mathcal{F}\left(\bw_{h}, \bu_{h}, \bH_h, \bv_{h}\right) = a\left(\bw_{h}, \bu_{h}, \bv_{h}\right) + \Rn^{-1}\ip{\nabla\times\bw_{h}}{\bv_{h}}{\Omega}- \mathsf{c}\ a\left(\nabla\times\bH_{h},\bH_{h},\bv_{h}\right),\]
\[\mathcal{M}\left(\bu_{h}, \bH_{h},\boldsymbol{b}_{h}\right)=\Rm^{-1}\ip{\nabla\times \bH_{h}}{\nabla\times\boldsymbol{b}_{h}}{\Omega} -\ a\left(\bu_{h}, \bH_{h}, \nabla\times\boldsymbol{b}_{h}\right).\]
Repeating the analyses in Section~\ref{SUBSEC: semi-discrete conservation} will reveal that this fully discrete formulation strongly preserves conservation of mass and conservation of charge, and weakly satisfies Gauss's law for magnetism. Especially, compared to the decoupled formulation, it preserves the correct dissipation of energy, which leads to conservation of energy for the ideal MHD (for example, see the bottom-left diagram of Fig.~\ref{fig:Energy tests}).
}

\bibliographystyle{elsarticle-num}
\bibliography{ref}

\begin{thebibliography}{10}
\expandafter\ifx\csname url\endcsname\relax
  \def\url#1{\texttt{#1}}\fi
\expandafter\ifx\csname urlprefix\endcsname\relax\def\urlprefix{URL }\fi
\expandafter\ifx\csname href\endcsname\relax
  \def\href#1#2{#2} \def\path#1{#1}\fi

\bibitem{Davidson_2016}
P.~A. Davidson, Introduction to Magnetohydrodynamics, 2nd Edition, Cambridge
  Texts in Applied Mathematics, Cambridge University Press, 2016.

\bibitem{BENSALAH20015867}
N.~{Ben Salah}, A.~Soulaimani, W.~G. Habashi, A finite element method for
  magnetohydrodynamics, Computer Methods in Applied Mechanics and Engineering
  190~(43) (2001) 5867--5892.

\bibitem{ZHANG202345}
X.~Zhang, H.~Su, A decoupled, unconditionally energy-stable and
  structure-preserving finite element scheme for the incompressible {MHD}
  equations with magnetic-electric formulation, Computers \& Mathematics with
  Applications 146 (2023) 45--59.

\bibitem{hu2017divB}
K.~{Hu}, Y.~{Ma}, J.~{Xu}, Stable finite element methods preserving
  $\mathrm{div}\cdot \boldsymbol{B}=0$ exactly for {MHD} models, Numerische
  Mathematik 135 (2017) 371--396.

\bibitem{HU2021110284}
K.~Hu, Y.-J. Lee, J.~Xu, Helicity-conservative finite element discretization
  for incompressible {MHD} systems, Journal of Computational Physics 436 (2021)
  110284.

\bibitem{BRACKBILL1980426}
J.~Brackbill, D.~Barnes, The effect of nonzero $\nabla\cdot\boldsymbol{B}$ on
  the numerical solution of the magnetohydrodynamic equations, Journal of
  Computational Physics 35~(3) (1980) 426--430.

\bibitem{DING2024}
Q.~Ding, X.~Long, S.~Mao, R.~Xi, Second order unconditionally convergent fully
  discrete scheme for incompressible vector potential mhd system, Journal of
  Scientific Computing 100 (2024).

\bibitem{TOTH2000605}
G.~Tóth, The $\nabla\cdot\boldsymbol{B}=0$ constraint in shock-capturing
  magnetohydrodynamics codes, Journal of Computational Physics 161~(2) (2000)
  605--652.

\bibitem{Balsara_2004}
D.~S. Balsara, J.~Kim, A comparison between divergence-cleaning and
  staggered-mesh formulations for numerical magnetohydrodynamics, The
  Astrophysical Journal 602~(2) (2004) 1079.

\bibitem{Evans1988SimulationOM}
C.~R. Evans, J.~F. Hawley, Simulation of magnetohydrodynamic flows: A
  constrained transport method, The Astrophysical Journal 332 (1988) 659--677.

\bibitem{refId0}
X.~Zhang, H.~Su, X.~Li, A fully discrete finite element method for a
  constrained transport model of the incompressible {MHD} equations, ESAIM:
  M2AN 57~(5) (2023) 2907--2930.

\bibitem{Heister2017}
T.~Heister, M.~Mohebujjaman, L.~G. Rebholz, Decoupled, unconditionally stable,
  higher order discretizations for {MHD} flow simulation, Journal of Scientific
  Computing 71 (2017) 21–43.

\bibitem{ZHANG2022110752}
G.-D. Zhang, X.~He, X.~Yang, A fully decoupled linearized finite element method
  with second-order temporal accuracy and unconditional energy stability for
  incompressible {MHD} equations, Journal of Computational Physics 448 (2022)
  110752.

\bibitem{10.1093/imanum/drad005}
J.~Droniou, L.~Yemm, {A hybrid high-order scheme for the stationary,
  incompressible magnetohydrodynamics equations}, IMA Journal of Numerical
  Analysis 44~(1) (2023) 262--296.

\bibitem{GAWLIK2022110847}
E.~S. Gawlik, F.~Gay-Balmaz, A finite element method for {MHD} that preserves
  energy, cross-helicity, magnetic helicity, incompressibility, and
  $\mathrm{div}\cdot \boldsymbol{B}=0$, Journal of Computational Physics 450
  (2022) 110847.

\bibitem{LAAKMANN2023112410}
F.~Laakmann, K.~Hu, P.~E. Farrell, Structure-preserving and helicity-conserving
  finite element approximations and preconditioning for the {Hall} {MHD}
  equations, Journal of Computational Physics 492 (2023) 112410.

\bibitem{MA202228}
H.~Ma, P.~Huang, A vector penalty-projection approach for the time-dependent
  incompressible magnetohydrodynamics flows, Computers \& Mathematics With
  Applications 120 (2022) 28--44.

\bibitem{ZHANG2024210}
Y.~Zhang, X.~Feng, H.~Su, Fully decoupled, linear and unconditionally energy
  stable time discretization scheme for solving the unsteady thermally coupled
  magnetohydrodynamic equations with variable density, Applied Numerical
  Mathematics 197 (2024) 210--229.

\bibitem{ZHANG2024113080}
Y.~Zhang, A.~Palha, M.~Gerritsma, Q.~Yao, A {MEEVC} discretization for
  two-dimensional incompressible {Navier-Stokes} equations with general
  boundary conditions, Journal of Computational Physics 510 (2024) 113080.

\bibitem{Oden2010}
J.~T. Oden, L.~F. Demkowicz, Applied Functional Analysis, Second Edition,
  Taylor \& Francis, 2010.

\bibitem{Arnold_Falk_Winther_2006}
D.~N. Arnold, R.~S. Falk, R.~Winther, Finite element exterior calculus,
  homological techniques, and applications, Acta Numerica 15 (2006) 1–155.

\bibitem{boffi2013mixed}
D.~Boffi, F.~Brezzi, M.~Fortin, Mixed Finite Element Methods and Applications,
  Springer, 2013.

\bibitem{buffa2001traces}
A.~Buffa, P.~Ciarlet~Jr., On traces for functional spaces related to
  {M}axwell's equations {P}art {I}: {A}n integration by parts formula in
  {L}ipschitz polyhedra, Mathematical Methods in the Applied Sciences 24~(1)
  (2001) 9--30.
\newblock \href
  {https://doi.org/10.1002/1099-1476(20010110)24:1<9::AID-MMA191>3.0.CO;2-2}
  {\path{doi:10.1002/1099-1476(20010110)24:1<9::AID-MMA191>3.0.CO;2-2}}.

\bibitem{dao_structure_2024}
T.~A. Dao, M.~Nazarov, I.~Tomas, A structure preserving numerical method for
  the ideal compressible {MHD} system, Journal of Computational Physics 508
  (2024) 113009.
\newblock \href {https://doi.org/10.1016/j.jcp.2024.113009}
  {\path{doi:10.1016/j.jcp.2024.113009}}.

\bibitem{wu_provably_2018}
K.~Wu, C.-W. Shu, A {Provably} {Positive} {Discontinuous} {Galerkin} {Method}
  for {Multidimensional} {Ideal} {Magnetohydrodynamics}, SIAM Journal on
  Scientific Computing 40~(5) (2018) B1302--B1329, publisher: Society for
  Industrial and Applied Mathematics.
\newblock \href {https://doi.org/10.1137/18M1168042}
  {\path{doi:10.1137/18M1168042}}.

\bibitem{brackbill_effect_1980}
J.~U. Brackbill, D.~C. Barnes, The {Effect} of {Nonzero} $\nabla\cdot$ {B} on
  the numerical solution of the magnetohydrodynamic equations, Journal of
  Computational Physics 35~(3) (1980) 426--430.
\newblock \href {https://doi.org/10.1016/0021-9991(80)90079-0}
  {\path{doi:10.1016/0021-9991(80)90079-0}}.

\bibitem{ramshaw_method_1983}
J.~D. Ramshaw, A method for enforcing the solenoidal condition on magnetic
  field in numerical calculations, Journal of Computational Physics 52~(3)
  (1983) 592--596.
\newblock \href {https://doi.org/10.1016/0021-9991(83)90009-8}
  {\path{doi:10.1016/0021-9991(83)90009-8}}.

\bibitem{balsara_staggered_1999}
D.~Balsara, D.~Spicer, A {Staggered} {Mesh} {Algorithm} {Using} {High} {Order}
  {Godunov} {Fluxes} to {Ensure} {Solenoidal} {Magnetic} {Fields} in
  {Magnetohydrodynamic} {Simulations}, Journal of Computational Physics 149~(2)
  (1999) 270--292.
\newblock \href {https://doi.org/10.1006/jcph.1998.6153}
  {\path{doi:10.1006/jcph.1998.6153}}.

\bibitem{ZHANG2022110868}
Y.~Zhang, A.~Palha, M.~Gerritsma, L.~G. Rebholz, A mass-, kinetic energy- and
  helicity-conserving mimetic dual-field discretization for three-dimensional
  incompressible {Navier-Stokes} equations, {part I: Periodic domains}, Journal
  of Computational Physics 451 (2022) 110868.

\bibitem{Hairer_Lubich_Wanner_2003}
E.~Hairer, C.~Lubich, G.~Wanner, Geometric numerical integration illustrated by
  the {Störmer–Verlet} method, Acta Numerica 12 (2003) 399–450.

\bibitem{hairer2006geometric}
E.~Hairer, C.~Lubich, G.~Wanner, Geometric numerical integration:
  Structure-preserving algorithms for ordinary differential equations, Vol.~31,
  Springer Science \& Business Media, 2006.

\bibitem{shen2018sav}
J.~Shen, J.~Xu, J.~Yang, The scalar auxiliary variable ({SAV}) approach for
  gradient flows, Journal of Computational Physics 353 (2018) 407--416.
\newblock \href {https://doi.org/10.1016/j.jcp.2017.10.021}
  {\path{doi:10.1016/j.jcp.2017.10.021}}.

\bibitem{Kreeft2011}
J.~Kreeft, A.~Palha, M.~Gerritsma, {Mimetic framework on curvilinear
  quadrilaterals of arbitrary order}, arXiv:1111.4304 (2011) 69.

\bibitem{zhang2022phd}
Y.~Zhang, {Mimetic Spectral Element Method and Extensions toward Higher
  Computational Efficiency} (2022).

\bibitem{nedelec}
J.~C. N{\'{e}}d{\'{e}}lec, Mixed finite elements in {$\mathbb{R}^3$}, Numer.
  Math. 35 (1980) 315--341.

\bibitem{raviart606mixed}
P.~A. Raviart, J.~M. Thomas, {A mixed finite element method for 2nd order
  elliptic problems}, Mathematical Aspects of the Finite Element Method,
  Lecture Notes in Mathematics 606 (1977) 292--315.

\bibitem{carstensen2016breaking}
C.~{Carstensen}, L.~{Demkowicz}, J.~{Gopalakrishnan}, Breaking spaces and forms
  for the {DPG} method and applications including {M}axwell equations,
  Computers \& Mathematics With Applications 72~(3) (2016) 494--522.

\bibitem{doi:10.1137/1.9781611977738}
L.~F. Demkowicz, Mathematical Theory of Finite Elements, Society for Industrial
  and Applied Mathematics, Philadelphia, PA, 2023.

\bibitem{Orszag_Tang_1979}
S.~A. Orszag, C.-M. Tang, Small-scale structure of two-dimensional
  magnetohydrodynamic turbulence, Journal of Fluid Mechanics 90~(1) (1979)
  129–143.

\bibitem{kraus2018variational}
M.~Kraus, O.~Maj, Variational integrators for ideal magnetohydrodynamics
  (2018).

\end{thebibliography}

%
%
%
%

\end{document}